\documentclass[smallextended,final]{svjour3}   

\usepackage{mathptmx}       
\usepackage{helvet}         
\usepackage{courier}        
\usepackage{makeidx}         
\usepackage{graphicx}        
\usepackage{multicol}        
\usepackage[bottom]{footmisc}

\makeindex             

\usepackage{tocloft}
\cftpagenumbersoff{figure}
\cftpagenumbersoff{table}

\usepackage{marvosym}

\usepackage{caption}

\usepackage{amsfonts,amssymb,amsmath,textcomp}
\usepackage{color}
\usepackage{listings}
\usepackage{xspace}
\usepackage[utf8x]{inputenc}
\usepackage[english]{babel}
\usepackage[active]{srcltx}

\usepackage{natbib}
\usepackage[title,titletoc]{appendix}

\usepackage[backref=true,bookmarksnumbered=true,bookmarks,colorlinks=true]{hyperref}
\hypersetup{urlcolor=blue, citecolor=blue,linkcolor=blue}

\newcommand{\DIV}{\mathop{\mathrm{div}}\nolimits}
\newcommand{\rot}{\mathop{\mathrm{curl}}\nolimits}
\newcommand{\grad}{\mathop{\mathrm{grad}}\nolimits}
\newcommand{\op}[2]{\mathop{\mathrm{{#1}_{#2}}}\nolimits}
\newcommand{\prt}[2]{\frac{\partial{#1}}{\partial{#2}}}
\newcommand{\prtt}[2]{\frac{\partial^2{#1}}{\partial{#2^2}}}

\newcommand{\prtq}[3]{\frac{\partial^2{#1}}{\partial{#2}\partial{#3}}}

\newcommand{\vE}{{\mathbf{E}}}
\newcommand{\vtE}{\mathbf{\widetilde E}}
\newcommand{\vH}{{\mathbf H}}

\newcommand{\vJ}{{\mathbf J}}
\newcommand{\vU}{{\mathbf U}}
\newcommand{\vV}{{\mathbf V}}
\newcommand{\vW}{{\mathbf W}}
\newcommand{\AD}{\mathrm{\Omega}}

\newcommand{\sW}{{\widetilde{\mathbf W}}^N}

\newcommand{\sL}[1]{{\mathbf L_2}[#1]}
\newcommand{\GE}{\widehat {G}^E}
\newcommand{\GH}{\widehat {G}^H}
\newcommand{\opp}[3]{\mathop{\mathrm{{\mathbf{#1}_{#2}^{#3}}}}\nolimits{}}
\newcommand{\Gem}{\opp{G}{E}{m}}
\newcommand{\PN}{\mathop{\mathrm{{\mathbf{P}}^{N}}}\nolimits}

\newcommand{\vAj}[2]{\mathbf{ #1}_{{#2}}}

\newcommand{\prtqq}[3]{\frac{\partial^{#1}}{\partial{#2}\partial{#3}}}
\newcommand{\tsig}{{\sigma}}
\newcommand{\conj}[1]{\overline{#1}}
\newcommand{\sqresig}{\sqrt{\Re{\tsig_b}}}
\renewcommand{\Re}[1]{\mathop{\mathbf{Re}}{#1}}

\newcommand{\MKCOMMENT}[1]{}
\begin{document}

\title{High-performance Parallel Solver for Integral Equations of Electromagnetics  Based on Galerkin Method 
}

\titlerunning{High-performance Parallel IE Solver Based on Galerkin Method}        

\author{Mikhail Kruglyakov       \and
       Lidia Bloshanskaya 
}


\institute{Mikhail Kruglyakov (\Letter)
	  \at Faculty of Computational Mathematics and Cybernetics, Lomonosov MSU, GSP-1,\\ Leninskiye Gory, 1-52, 119991 Moscow, Russia \\Tel. +74959391919  Fax. +74959392596
	  \at Institute of Geophysics, ETH Zurich, Sonneggstrasse 5, 8092 Zurich, Switzerland \\
	  \email{mkruglyakov@cs.msu.su}
	  \and 
	  Lidia Bloshanskaya 
	    \at SUNY New Paltz, Department of Mathematics, 1 Hawk Dr, New Paltz, NY  12561,  U.S.A. \\
	    \email{bloshanl@newpaltz.edu}
  }

  \date{Received: \hphantom{\today ~}/ Accepted: \hphantom{\today ~}}

\maketitle

\begin{abstract} 
	A new parallel  solver for the  volumetric integral equations (IE)  of electrodynamics is presented. The solver is based on the Galerkin method which ensures the convergent numerical solution. The main features  include: (i) the memory usage is  8 times lower, compared to analogous IE based algorithms, without additional restriction on the background media; (ii) accurate and stable method to compute  matrix coefficients corresponding to the IE; (iii)  high degree of parallelism. The solver's computational efficiency is shown on a  problem of magnetotelluric sounding  of the high conductivity contrast media.  A good agreement with the results obtained with the second order finite element method is demonstrated. Due to effective approach to parallelization and distributed data storage the program exhibits perfect scalability on different hardware platforms.

 \keywords{Integral equations \and Forward modeling \and Electromagnetic sounding \and Galerkin method \and Green's tensor \and High-performance computing}
\end{abstract}

\section{Introduction} 

Electromagnetic (EM) methods of geophysics are used to model the subsurface   electrical conductivity distribution. 
Conductivity is affected by the rock type and composition, temperature, and fluid/melt content and thus can be used in various engineering and industrial problems like  detecting hydrocarbon (low-conductive) and geothermal or ore (high-conductive) reservoirs.
Measured  electrical and/or magnetic fields  are further interpreted via the calculations for a given three-dimensional  model of conductivity distribution. Maxwell's equation describing the EM field distribution can not be solved analytically in general  case requiring numerical simulation. Large number of such simulations is required, and complex large scale  models are invoked, \cite{ChaveJonesMT}. 

The  growing amount of data calls for the development of new  numerical methods capable to deliver fast and accurate EM simulations and harness the computational power provided by modern high-performance multi-core and multi-node platforms. 

There are three basic  approaches  to the numerical simulation of EM fields in the conductive media: finite-difference (FD), finite-element (FE) and volumetric integral equation (IE) methods.  The FD schemes dominated in EM sounding for decades, \cite{Mackie1994, Haber2001, Newman2001, egbert2012,Jaysaval01122015}. However,  the FE methods  has become popular in recent years, \cite{Schwarzbach2011, Farquharson2011, Puzyrev2013, Ren2013, Grayver2015}. 
The IE methods are not so common due to the certain  difficulties in their implementation. For their usage one can refer to \cite{avdeev2002, Hursan02, sing08, koya08, Kamm2014}.

The main difference in these approaches is in the discretization outcome. The model usually consists of a number of the non uniform three-dimensional  anomalies embedded in the one-dimensional (layered) background media.
While the FD and FE methods produce large sparse systems, \cite{Ernst2011}, the~IE method results in the compact dense system matrices.
The compactness is attained since the modeling region is confined only to the three-dimensional conductivity structures (anomalies) under the investigation, \cite{Raiche1974, Weidelt1975}.
Note that boundary conditions are satisfied by construction of the Green's functions. By contrast, in the FD and FE methods one has to discretize a volume much bigger in both  lateral and vertical directions in order to enable the decay (or stabilization) of the EM field at the boundaries of the modeling domain, \cite{Grayver2015, Mulder2006}.  Another distinction between the methods is in the condition number of matrices (it controls the stability of the solution).  In FD and FE methods the condition number depends on the discretization and frequency, whereas in IE approach it does not, \cite{pankratov1995, sing95}.

The main focus in development and implementation of the numerical methods for the IE is the efficiency and performance for the large number of unknown parameters. 
The proposed new iterative numerical solver for the IE  addresses this issue both mathematically (increasing the accuracy and stability of coefficient computations and reducing the memory usage) and computationally (allowing high degree of parallelism without memory loss at the nodes).
The special class of integral equations is used:  the integral equations with contracting kernel (CIE), \cite{pankratov1995, sing95}.
The CIE were proved to have a unique solution and a well-conditioned
(by construction) system matrix, \cite{pankratov1995, sing95}. The Galerkin method is used to solve CIE numerically, \cite{Delves_Mohamed-Galerkin}. In \cite{kruglyakov2011, sing08}   it was proved that the corresponding solution converges. 

Two main issues in the IE numerical solution are the calculation of matrix coefficients with sufficient accuracy, \cite{Wannamaker1991}, and the storage of this dense matrix in some packed form, \cite{AvdeevEtAl1997,avdeev2002}. The Galerkin method with piece-wise constant basis allows to address both of these issues.  Namely, the system matrix is decomposed into sums and products of diagonal  and block-Toeplitz matrices. First, matrix coefficients (i.e., double volumetric integrals of the product of basis functions and CIE kernel) are analytically transformed into the one-dimensional convolution integrals. These convolution integrals are then computed by digital filtering approach. Keep in mind, that only weights of these filters are computed numerically, whereas the functions in the knots are computed analytically. The resulting system is solved using the Flexible GMRES, \cite{Saad1993}.
 
The proposed solver exhibits three main features.
1) The memory usage is  8 times lower, compared to the analogous IE solving algorithms. It is important to stress that the memory is saved for any background with an arbitrary number of layers. In~contrast, in  \cite{Kamm2014} the memory reduction is achieved only in case of  homogenous half-space as a background and for uniform vertical discretization. In \cite{Avdeev09,koya08,Sun2015} it is achieved at the expense of accuracy and performance. 
The idea behind the proposed solver is to combine the Galerkin method with the properties of EM field, namely Lorentz reciprocity. The matrix of the ensuing  linear system can be separated in symmetric and antisymmetric submatricies. This reduces the memory requirements by 8 times,  Sect.~\ref{sec_features}. 
2) Efficient and accurate method for the  computation of these matrices,  Sect.~\ref{sec_features}. 
3) The  implementation  with high degree of parallelism. The computational experiments performed with “Bluegene” and “Lomonosov” supercomputers from MSU, and  high-performance computer (HPC) ``Piz Daint'' from Swiss National Supercomputing Center  show that the solver makes the best usage of 128 to 2,048 nodes for calculation at a single frequency and a single source.  The program exhibits perfect scalability.

This paper is organized as follows. Section~\ref{ContractionIE} is devoted to the overview of the CIE  approach and the construction of the approximating system of linear equations. Section \ref{sec_features} addresses the reduction in memory requirements, the computation of matrix coefficients, and features of parallel implementation. 
In Sect.~\ref{Sec_Example} the computational results for high (more than $3 \cdot 10^4$) conductivity contrast  COMMEMI3D-3 model (\cite{Hursan02,Varentsov2000}) are compared with the corresponding results obtained using  FE method  by \cite{Grayver2015}.  
The Appendices~\ref{AppendixGreensTensor} to \ref{AppendixHorizontalIntegration} provide the mathematical details of the presented method.

\section{Contracting Integral Equation}\label{ContractionIE} 

\subsection{Overview}
Assume that  the EM fields are induced by the external electric currents~$\vJ_{ext}$. 
Moreover, assume that the EM fields are time dependent as $e^{-i\omega t}$, where $\omega$ is an angular frequency, $i=\sqrt{-1}$ and magnetic permeability $\mu_0$ is the same in whole space.
Let $\sigma(M)$, $\Re{\sigma(M)} \ge 0$ be a three-dimensional complex conductivity distribution in  space.  Then the electrical field $\vE$ and the magnetic field $\vH$ give the solution of the system of Maxwell's equations
\begin{equation} 
	\label{Maxwell}
	\left\{
		\begin{aligned}
			&\rot \vH =\tsig \vE+\vJ_{ext},\\
			&\rot \vE=i\omega\mu_0 \vH.\\
		\end{aligned}
	\right.
\end{equation}
The solution of~\eqref{Maxwell} is unique under the additional radiation conditions at infinity,  \cite{Hohmann1988}.

\begin{figure}[t!]
	\centering \includegraphics[width=\textwidth, keepaspectratio]{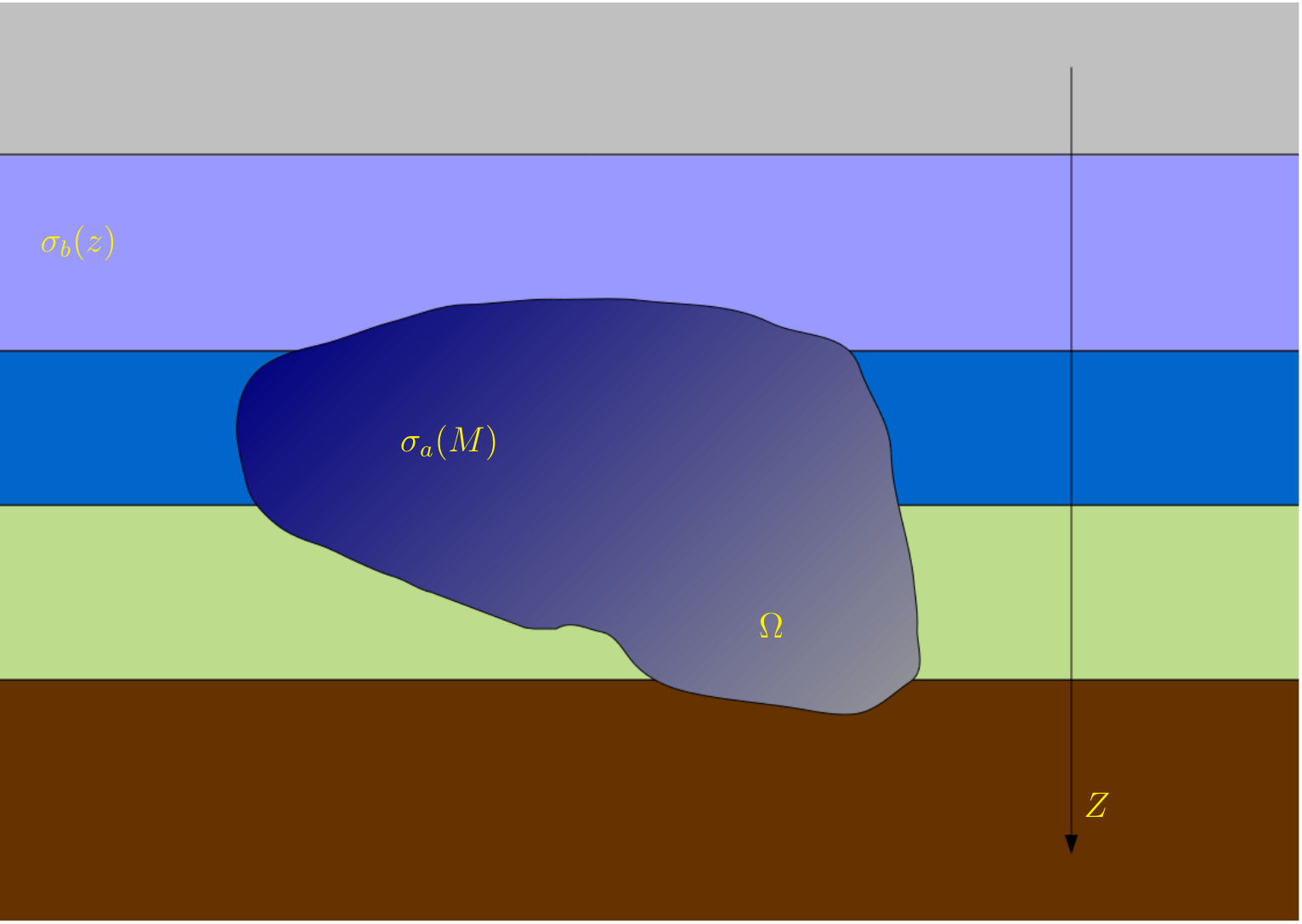}
\caption{Typical model}
\label{typical_model}
\end{figure}

Let $\AD\subset\mathbf{R}^3$ be some bounded domain  and $\tsig(M)=\tsig_b(z)$ for $M(x,y,z) \not \in \AD$ and $\tsig(M)=\tsig_a(M)$~for~$M \in \AD$~(Fig.~\ref{typical_model}).
Then for any $M \in \mathbf{R}^3$ the fields $\vE(M)$ and $\vH(M)$  are expressed in terms of the integrals
\begin{equation} 
	\label{cont_formulas}
	\begin{aligned}
		\vE(M)&=\vE^N(M)+\int\limits_{\AD}\GE(M,M_0)\left(\tsig_a(M_0)-\tsig_b(M_0)\right)\vE(M_0)d\AD_{M_0},\\
		\vH(M)&=\vH^N(M)+\int\limits_{\AD}\GH(M,M_0)\left(\tsig_a(M_0)-\tsig_b(M_0)\right)\vE(M_0)d\AD_{M_0}.\\
	\end{aligned}
\end{equation}
Here $\GE$, $\GH$ are electrical and magnetic Green's tensors respectively, \cite{dmitriev2002, pankratov1995}. The terms $\vE^N$, $\vH^N$ are called the normal electric and magnetic fields, corresponding. They form the  solution of the  system
\begin{equation} 
	\label{Maxwell2}
	\left\{
		\begin{aligned}
			&\rot \vH^N =\tsig_b(z) \vE^N+\vJ_{ext},\\
			&\rot \vE^N=i\omega\mu_0 \vH^N\\
		\end{aligned}
	\right.
\end{equation}
with corresponding conditions at the infinity. 
Note, that $\GE$, $\GH$ are indepent  of  the anomalous conductivity $\tsig_a$.

Let $\sL{\AD}$ be a Hilbert functional space of vector functions $\vV$ with the following norm and dot product
\begin{equation}
	\label{L_defininition}
	\begin{aligned}
		&\left(\vV,\vU\right)_{\sL{\AD}}=\int\limits_{\AD}\left(V_x(M)\conj{U}_x(M)+V_y(M)\conj{U}_y(M)+V_z(M)\conj{U}_z(M)\right)d\AD_M,\\
		&\left\|\vV\right\|_{\sL{\AD}}=\sqrt{\left(\vV,\vV\right)}.
	\end{aligned}
\end{equation}
Suppose $\Re{\tsig_b}(z)>0$ for $M(x,y,z) \in \AD$ (the typical EM sounding situation), then the operator $\Gem$ is defined as
\begin{equation}
 \label{Grmod}
 \Gem \vAj{V}{}=\sqrt{\Re{\tsig_b}}\op{\widehat G^E}{}\left[2\sqrt{\Re{\tsig_b}}\vAj{V}{}\right]+\vAj{V}{},
\end{equation} 
where $\op{\widehat G^E}{}$ is an integral operator from the first equation in \eqref{cont_formulas}. The operator $\Gem$ is a contracting operator in $\sL{\AD}$,  \cite{pankratov1995, sing95}. Using~\eqref{cont_formulas}~and~\eqref{Grmod} one obtains the CIE for $\vE$
\begin{equation}
\label{iEH320op}
\begin{aligned}
	&\left(\op{I}{}-\Gem\frac{b}{a}\right)\vtE=\sqrt{\Re{\tsig_b}}\vE^N,\\
	&\vtE=a\vE,~ ~ ~ a=\frac{\tsig_a+\conj{\tsig_b}}{2\sqresig}~~~~b=\frac{\tsig_a-\tsig_b}{2\sqresig},
\end{aligned}
\end{equation} 
where $\op{I}{}$ is the identity operator and $\conj{\tsig_b}$ means complex conjugation of ${\tsig_b}$.

\subsection{Galerkin Method} 
Suppose the domain $\AD$ is divided in nonoverlapping subdomains $\AD=\cup \AD_n$, $n=1\dots N$ and $\tsig_b(M)=\tsig_b^n$, $\tsig_a(M)=\tsig_a^n$ for $M \in \AD_n, n=1\dots N$. For each subdomain $\AD_n$  define the function $W_n(M)$ as
\begin{equation}
	\begin{aligned}
			W_n(M)=\left\{\begin{aligned}
					\frac{1}{V_n},\ M \in \AD_n,\\
					0,\ M \not \in \AD_n,
					\end{aligned}\right.&~ ~ ~
			V_n=\int\limits_{\AD_n}d\AD_M,& n=1\dots N.\\
	\end{aligned}
\end{equation}
Let $\sW$ be a linear span of the vector functions $\vW_n$, $\vW_n=(W_{n_x},W_{n_y},W_{n_z})$, $n_x,n_y,n_z=1\dots N$  and $\PN$ be a projection operator from $\sL{\AD}$ to $\sW$
\begin{equation}
	\label{p_def}
	\forall \vAj{F}{}\in \sL{\AD}~ ~ ~\left[\PN[\vAj{F}{}]\right]_\gamma=\sum\limits_{n=1}^N\alpha_n^\gamma W_n,~ ~ ~ \alpha_n^\gamma=\frac{\int\limits_{\AD_n}F_\gamma(M)d\AD_M}{V_n},
\end{equation}
where $\gamma={x,y,z}$. Note that $\left\| \PN \right\|=1$.

Applying $\PN$ to the first equation in \eqref{iEH320op}  one obtains the  operator equation in $\sW$
\begin{equation}
\label{pr_op}
\begin{aligned}
	&\vW-\PN\Gem\frac{b}{a}\vW=\vW^0,\\
	&\vW^0=\PN\sqresig\vE^N.\\
\end{aligned}
\end{equation}

Since  $\frac{b}{a}<1$, $\Gem$ is a contracting operator and $\left\|\PN\right\|=1$, it can be easily shown that~\eqref{pr_op} has a unique solution $\vW$ in $\sW$,  \cite{kruglyakov2011, sing08}. Moreover, $\vW$ approximates $\vE$ with the first order of $d$ in $\sL{\AD}$, where $d=\max\limits_{n=1\dots N} d_n$, $d_n$ is a diameter of the subdomain~$\AD_n$,~\cite{kruglyakov2011}.

Using the definition of $\sW$  the components of $\vW=(W_x,W_y,W_z)$ can be expressed as
\begin{equation}
	\label{w2u}
	W_\gamma=\sum_{k=1}^NU_n^\gamma W_n,~ ~ ~\gamma=x,y,z.
\end{equation}
Using \eqref{p_def} to \eqref{w2u}  and taking into account that $\tsig_a,\tsig_b$ are piecewise functions,  one obtains the following system of linear equations for the coefficients ${\vU_n=(U_n^x,U_n^y,U_n^z)}$, ${n=1\dots N}$
\begin{equation}
	\label{slae1}
	\vU_n-\sum_{m=1}^{N}\hat \gamma^m\widehat K_n^m\vU_m=\vU_n^0,
\end{equation}
where
\begin{equation}
	\label{slae_matricies}
	\begin{aligned}
		\widehat K_n^m&=\hat I+\frac{2}{V_n}\sqrt{\Re\tsig_b^m\Re\tsig_b^n}\hat B_n^m,\\
		\widehat B_n^m&=\int\limits_{\AD_n}\int\limits_{\AD_m}\GE(M,M_0)d\AD_{M_0}d\AD_{M}, \\
		\hat \gamma^m&=\frac{\tsig_a^m-\tsig_b^m}{\tsig_a^m+\conj{\tsig_b^m}},\\
		\vU_n^0&=\frac{\sqrt{\Re\tsig_b^n}}{V_n}\int\limits_{\AD_n} \vE^N(M)d\AD_M.\\
	\end{aligned}
\end{equation}
Note that $\widehat K_n^m$, $\widehat B_n^m$, $\hat I$, $\hat \gamma^m$ are  $3 \times 3$ matrices, $\hat I$ is an identity matrix, $\hat \gamma^m$ is a diagonal matrix. 
The system \eqref{slae1} has a unique solution, \cite{kruglyakov2011, sing08}.

Using the solution $\vU_n$ of system \eqref{slae1} one can approximate $\widetilde \vE(M)$ and $\widetilde \vH(M)$ for any point $M \in R^3$ with
\begin{equation}
	\label{APR_lab}
\begin{aligned}
\vE(M) \approx \widetilde \vE(M)&=\vE^N(M)+\sum\limits_{n=1}^N\left(\tsig_a^n-\tsig_b^n\right)\vU_n\int\limits_{\AD_n}\GE(M,M_0)d\AD_{M_0},\\
\vH(M) \approx \widetilde \vH(M)&=\vE^H(M)+\sum\limits_{n=1}^N\left(\tsig_a^n-\tsig_b^n\right)\vU_n\int\limits_{\AD_n}\GH(M,M_0)d\AD_{M_0}.\\
\end{aligned}
\end{equation}
Relations \eqref{APR_lab} are first order approximations of $d$ in $\mathbf{C}$ for $M \not \in \AD$, ~\cite{kruglyakov2011}, and result in fast and relatively simple computations. The main challenge is to calculate matrix coefficients and  solve the system~\eqref{slae1}.

\section{Computational Challenges}\label{sec_features} 

\subsection{Memory Requirements} 
The main challenge of the integral equation approach is in  solving of the system of linear equations with dense matrices \eqref{slae1}. The storage of these  matrices in  RAM is also problematic. The standard approach, \cite{AvdeevEtAl1997}, is to use the property
	\begin{equation}  
		\label{convolv}
		\GE(M,M_0)=\GE(x-x_0,y-y_0,z,z_0).
	\end{equation}
For the implementation purposes consider now $\AD\subset R^3$ to be a rectangular domain. As before~$\AD$~is~divided in $N=N_xN_yN_z$ rectangular subdomains~$\AD_n$,~$n=1,\ldots,N$,  where  $N_x,N_y,N_z$ are the number of subdomains in $X,Y,Z$ directions respectively. Suppose also that each $\AD_n$ has the same size $h_x\times h_y$ in $XY$ plane. Then
	\begin{equation}  
		\label{bcyl}
		\widehat B_n^m=	\int\limits_{\AD_n}\int\limits_{\AD_m}\GE(M,M_0)d\AD_{M_0}d\AD_M=\widehat B_n^m\left(I^n_x-I^m_x,I^n_y-I^m_y,I^n_z,I^m_z\right),
	\end{equation}
	where $I^n_x,I^m_x \in \{1,2, \dots N_x\}$, $I^n_y,I^m_y\in \{1,2, \dots N_y\}$, $I^n_z,I^m_z \in \{1,2, \dots N_z\}$, $n,m=1\dots N$. Therefore  $\widehat B_n^m$ is a block Toeplitz matrix induced by the block vector $(C^y_{-(N_y-1)},$ $C^y_{-(N_y-2)}$, $\dots, C^y_{N_y-2}$, $C^y_{N_y-1})$. 
Each block $C^y_i$, $i=-(N_y-1)\dots N_y-1$ is also a block Toeplitz matrix and is induced by the block vector $(D^i_{-(N_x-1)},$ $D^i_{-(N_x-2)},\dots D^i_{N_x-2},D^i_{N_x-1})$. The $D^i_j$ is a $3 \times 3$ block matrix with the structure
	\begin{equation} 
		\label{Q_def}
		D^i_j=Q(i,j)=\begin{pmatrix}
			Q_{xx}&Q_{xy}&Q_{xz}\\
			Q_{yx}&Q_{yy}&Q_{yz}\\
			Q_{zx}&Q_{zy}&Q_{zz}
		\end{pmatrix}.
	\end{equation}
	Here $Q_{\alpha\beta}$ are the matrices of the order $N_z$, $\alpha,\beta={x,y,z}$, $i=-(N_y-1)\dots N_y-1$, $j=-(N_x-1)\dots N_x-1$. 
	
	Let $A$ be a matrix corresponding to the system of linear equations \eqref{slae1}. Then
	\begin{equation} 
		\label{matrix_convolv}
		A=S+R_1BR_2,
	\end{equation}
	where $S$, $R_1$, $R_2$ are the diagonal matrices; $B=\left\{\widehat B_n^m\right\}$ is the block Toeplitz matrix described above.

	In view of \eqref{matrix_convolv} it follows that only  $36\cdot N_xN_yN_z^2\cdot 16+O(N_xN_yN_z)$ bytes are required to store matrix $A$ in double precision. Using the equivalence $G_{xy}^E=G_{yx}^E$ this requirement can be reduced to $32\cdot N_xN_yN_z^2\cdot 16+ O(N_xN_yN_z)$ bytes as in \cite{AvdeevEtAl1997}.
	This memory requirement can be reduced in 8 times by  virtue of the following Lemmas.
	\begin{lemma}\label{lemm1}If $ \GE(M,M_0)$ is an electrical Green's tensor of any layered media, then  it possesses symmetric and antisymmetric properties in Cartesian coordinates along the vertical dimension
\begin{equation}
	\label{tensor_symm}
	\begin{aligned}
		G^E_{xx}(x-x_0,y-y_0,z,z_0)&=G^E_{xx}(x-x_0,y-y_0,z_0,z),\\ 
		G^E_{yy}(x-x_0,y-y_0,z,z_0)&=G^E_{yy}(x-x_0,y-y_0,z_0,z),\\
		G^E_{zz}(x-x_0,y-y_0,z,z_0)&=G^E_{zz}(x-x_0,y-y_0,z_0,z),\\
		G^E_{xy}(x-x_0,y-y_0,z,z_0)&=G^E_{yx}(x-x_0,y-y_0,z,z_0),\\
		G^E_{xy}(x-x_0,y-y_0,z,z_0)&=G^E_{xy}(x-x_0,y-y_0,z_0,z),\\
		G^E_{zx}(x-x_0,y-y_0,z,z_0)&=-G^E_{xz}(x-x_0,y-y_0,z_0,z),\\ 
		G^E_{zy}(x-x_0,y-y_0,z,z_0)&=-G^E_{yz}(x-x_0,y-y_0,z_0,z).\\ 
	\end{aligned}
\end{equation}
	\end{lemma}
	\begin{lemma}\label{lemm2}If $ \GE(M,M_0)$ is an electrical Green's tensor of any layered media, then  it possesses symmetric and antisymmetric properties in Cartesian coordinates along the lateral dimensions
	\begin{equation}
		\label{tensor_symm_2}
	\begin{aligned}
		G^E_{\alpha\alpha}(x-x_0,y-y_0,z,z_0)&=G^E_{\alpha\alpha}(x_0-x,y-y_0,z,z_0)=\\ 
		G^E_{\alpha\alpha}(x-x_0,y_0-y,z,z_0)&=G^E_{\alpha\alpha}(x_0-x,y_0-y,z,z_0),\\ 
		G^E_{xy}(x-x_0,y-y_0,z,z_0)&=-G^E_{xy}(x_0-x,y-y_0,z,z_0)=\\
		-G^E_{xy}(x-x_0,y_0-y,z,z_0)&=G^E_{xy}(x_0-x,y_0-y,z,z_0),\\
		G^E_{zx}(x-x_0,y-y_0,z,z_0)&=-G^E_{zx}(x_0-x,y-y_0,z,z_0)=\\ 
		G^E_{zx}(x-x_0,y_0-y,z,z_0)&=-G^E_{zx}(x_0-x,y_0-y,z,z_0),\\ 
		G^E_{zy}(x-x_0,y-y_0,z,z_0)&=G^E_{zy}(x_0-x,y-y_0,z,z_0)=\\ 
		-G^E_{zy}(x-x_0,y_0-y,z,z_0)&=-G^E_{zy}(x_0-x,y_0-y,z,z_0),\\ 
		\alpha&\in \{x,y,z\}.
	\end{aligned}
	\end{equation}
\end{lemma}
These Lemmas are trivial corollaries from Lorentz reciprocity, \cite{Hohmann1988} and   formulas for the Green's tensor components (See Appendix \ref{AppendixGreensTensor}).
	Relations \eqref{bcyl} and \eqref{tensor_symm} give
\begin{equation} 
	\label{matrix_symm}
	\begin{aligned}
			Q_{zx}&=-Q_{xz}^T, Q_{zy}=-Q_{yz}^T,\\
			Q_{xx}&=Q_{xx}^T, Q_{yy}=Q_{yy}^T, Q_{zz}=Q_{zz}^T,\\
			Q_{xy}&=Q_{xy}^T=Q_{yx}=Q_{yx}^T,\\
	\end{aligned}
 \end{equation}
where $^T$ indicates a matrix transpose.

Therefore, one needs to store only $Q_{xz},Q_{yz}$ and upper diagonal parts of $Q_{xx}$, $Q_{xy}$, $Q_{yy}$, $Q_{zz}$.
Moreover the values $Q(i,j)$ can be stored only for $i=0,\dots,N_y-1, j=0,\dots,N_x-1$, since \eqref{tensor_symm_2} allows to obtain these values for negative $i$ or $j$ from suitable symmetric/antisymmetric properties.

Thus only  $2\cdot N_xN_y N_z\cdot(2N_z+1) \cdot 16$  bytes are required to store $\widehat B_n^m$ which is 8 times less than the memory requirements in \cite{AvdeevEtAl1997,avdeev2002, Hursan02}. It is worth to stress again, that this is valid for any background layered media and without the conditions on the subdomains to be of the same vertical sizes.  \MKCOMMENT{something about uniform in lateral? I don't know}

\subsection{Matrix Coefficients Computation}  

The next computational challenge of the Galerkin approach is  the evaluation of the coefficients  $\widehat B_n^m$, $n,m =1\dots N$, that is the double volumetric integrals of the $\GE$ in the RHS of \eqref{bcyl}, with desired accuracy. The components of $\GE$ are the improper integrals containing the Bessel functions, Appendix~\ref{AppendixGreensTensor}. The integration in vertical direction is performed analytically using the fundamental function of layered media approach from \cite{dmitriev2002}, Appendix~\ref{AppendixVerticalIntegration}. The main problem, however, is the integration over the horizontal domains. In this case one  needs to compute the fifth-order integrals over the fast-oscillating functions.

 The integrals in \eqref{bcyl} are double volumetric ones, thus they have only weak singularity. Therefore, one can change  the order of integration and make an appropriate substitution and  convert the fifth-order integral  to a convolution with the specific kernel. Following the standard approach of convolution calculation the spectrum of this kernel is computed and the digital filter is constructed, Appendix~\ref{AppendixHorizontalIntegration}.
 
 It is important to emphasize that both the knots and the weights in the obtained filter significantly depend on the integration domains. On the contrary, the integration over different horizontal domains is completely data independent. This is used in parallel algorithm.
 The computational experiments demonstrate (Sect.~\ref{Sec_Example}) that the used filters provide suitable accuracy even  for the models with  high conductivity contrast.
 
\subsection{Parallel Implementation} 
The most essential part of any iterative method for  solving a system of linear equations is the matrix-vector multiplication. Since matrix $B$ is a block Toeplitz matrix, one can use the two-dimensional Fast Fourier Transform (FFT) to speed up this operation, \cite{AvdeevEtAl1997}.
Therefore, instead of matrices $Q(i,j)$ the  discrete Fourier  transformations $\widetilde Q(i,j)$ are stored. This requires the same amount of memory since the discrete Fourier transform preserves the symmetric/antisymmetric properties  of data.

The multiplication of block Toeplitz matrix $B$ on some vector $\vV \in \sW$ is performed via the following three-step algorithm:
\begin{enumerate}
	\item Compute $3 N_z$ forward   FFT of vector $\vV$;
	\item Compute $36N_xN_y$ algebraic matrix-vector multiplications of order $N_z$  to obtain vector~$\tilde \vV$;
	\item Compute $3 N_z$ backward  FFT of vector $\tilde \vV$.
\end{enumerate}
The multiplications in Step 2 are further divided into $4N_xN_y$ groups which are mutually data independent. This allows to implement the special scheme of distributed data storage and a parallel algorithm of IE solver, described below. 


For simplicity, consider $2N_y$ nodes and assume that $N_x$ is even. The distributed storage of matrix is organized in a special way: the half of block-vector $\widetilde Q(n,j)$, $j=0\dots N_x/2-1$ is stored at $n$th~node, {$n=0\dots N_y-1$}, while  $\widetilde Q(n,j)$, $j=N_x/2\dots N_x$ is stored at node $n+N_y$, $n=0\dots N_y-1$, Table \ref{mso}.
\begin{table}
\begin{center}
	\begin{tabular}{|c||c||c||c||c||c|}\hline
		Node 0		    &  $\dots$ & Node $N_y-1$& Node $N_y$ & \dots     & Node $2N_y-1$ \\\hline 
\vphantom{$\left\{\widetilde Q\right\}$}$ \widetilde Q(0,0)$ &  $\dots$ & $\widetilde Q(N_y-1,0) $ &$\widetilde Q(0,N_x/2) $& \dots& $\widetilde Q(N_y-1,N_x/2)$  \\

\vphantom{$\left\{\widetilde Q\right\}$}$ \widetilde Q(0,1)$  & \dots & $\widetilde Q(N_y-1,1) $ &$\widetilde Q(0,N_x/2+1) $& \dots& $\widetilde Q(N_y-1,N_x/2+1)$  \\

\vphantom{$\left\{\widetilde Q\right\}$}\vdots   & $\dots$ & \vdots  & \vdots & \dots& \vdots  \\

\vphantom{$\left\{\widetilde Q\right\}$}$ \widetilde Q(0,N_x/2-1)$ &  $\dots$ & $\widetilde Q(N_y-1,N_x-1) $ &$\widetilde Q(0,N_x-1) $& \dots& $\widetilde Q(N_y-1,N_x-1)$  \\\hline
	\end{tabular}
\end{center}
\caption{Matrix storage organization}
\label{mso}
\end{table}

\noindent This storage organization is used to develop the solver with suitable features of parallelization:
\begin{enumerate}
\item[(i)] The coefficients of matrices $B$, $S$, $R_1$, $R_2$ stored at different nodes are computed simultaneously and completely data independent;
\item[(ii)] The iterative method is executed using  the authors' distributed implementation of FGMRES \cite{Saad1993}, inspired by~\cite{Fraysse_et_al_2003};
\item[(iii)] The distributed two-dimensional Fourier transform is computed via the authors' implementation using FFTW3 library~\cite{FFTW05} for local FFT;
\item[(iv)] The calculation of the local algebraic matrix-vector multiplication  is processed by using  OpenBLAS library~\cite{OpenBlas} ; 
\item[(v)] For all the stages of the computational process the hybrid MPI+OpenMP scheme is used.
\end{enumerate}

\begin{figure}[t!]
	\begin{center}
		\includegraphics[width=\textwidth,keepaspectratio]{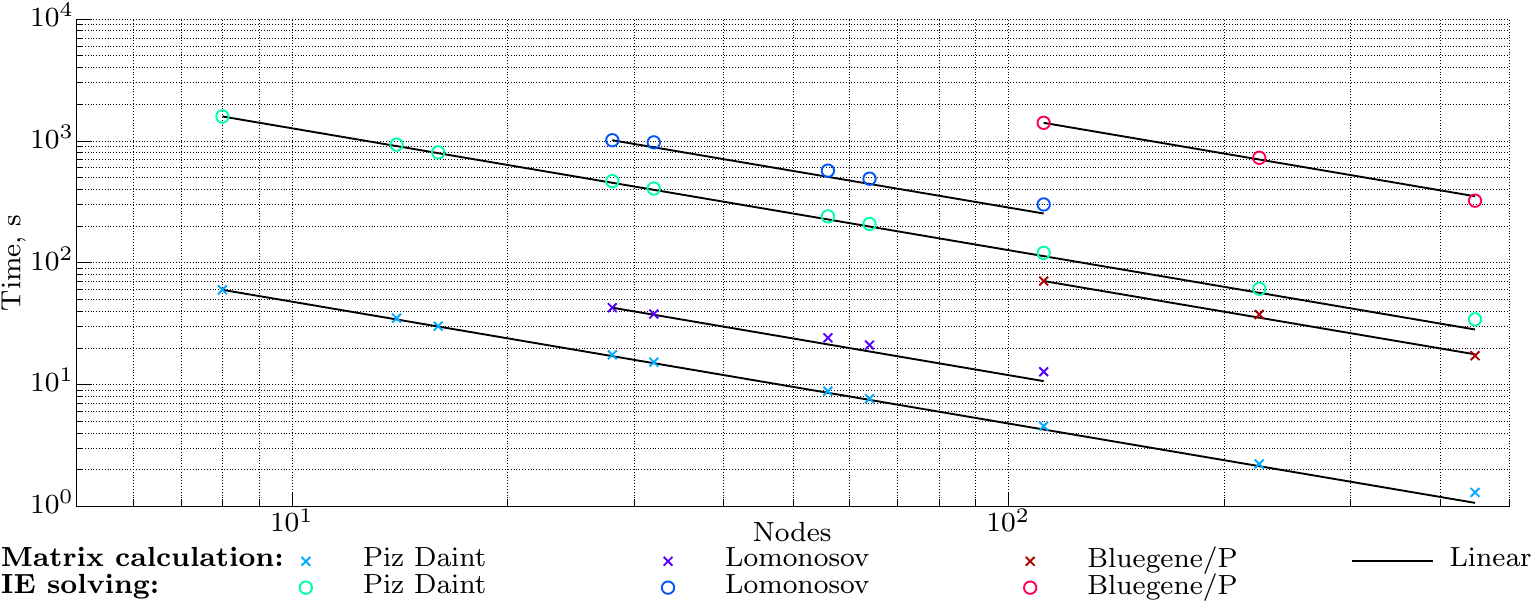}
\caption{Strong scalability for COMMEMI3D-3 model}
\label{par_time}
\end{center}
\end{figure}

 To demonstrate the scalability of the implemented parallelization the \newline COMMEMI3D-3  model is used  with $N_x=176$, $N_y=224$, $N_z=118$, that is with  cubic subdomains with 25\,m edges, Sect.~\ref{Sec_Example}. The computational experiments performed at ``Bluegene/P'', HPC ``Lomonosov'' (MSU) and Piz Daint (Swiss National Supercomputing Center)  showed good speed increment  depending on the number of processes~(Fig.~\ref{par_time}). Matrix calculation time includes time of FFT calculation of $\widehat B_n^m$. The solid black line means ideal linear speed up. Note, that for such high-contrast model matrix calculation time (crosses) is small enough compared to  solving of the system of linear equations (circles). One can see that the scalability is  close to a linear. 

\section{High Conductivity Contrast Modeling}
\label{Sec_Example}

The accurate computation of the EM field in a high conductivity contrast media is one of the most complex problems of EM modeling due to strong codependency between conductivity contrast and matrix condition number,  \cite{pankratov1995, Pankratov_Kuvshinov_2015, sing95}. The conductivity contrast means the ratio between the real parts of anomalous conductivity $\Re{\tsig_a(M)}$ and  background conductivity $\Re{\tsig_b(M)}$ at the same point $M$.

The high conductivity contrast  COMMEMI3D-3 model, \cite{Hursan02,Varentsov2000}, is used as one of the test models for the presented solver.  This model schematically describes the conductivity distribution typical for  the ore exploration by the audio-magnetotelluric sounding.  Following magnetotelluric (i.e., low-frequency) sounding tradition in the rest of this section the conductivity is a real-valued function.

\begin{figure}[t!]
	\centering \includegraphics[width=\textwidth, keepaspectratio]{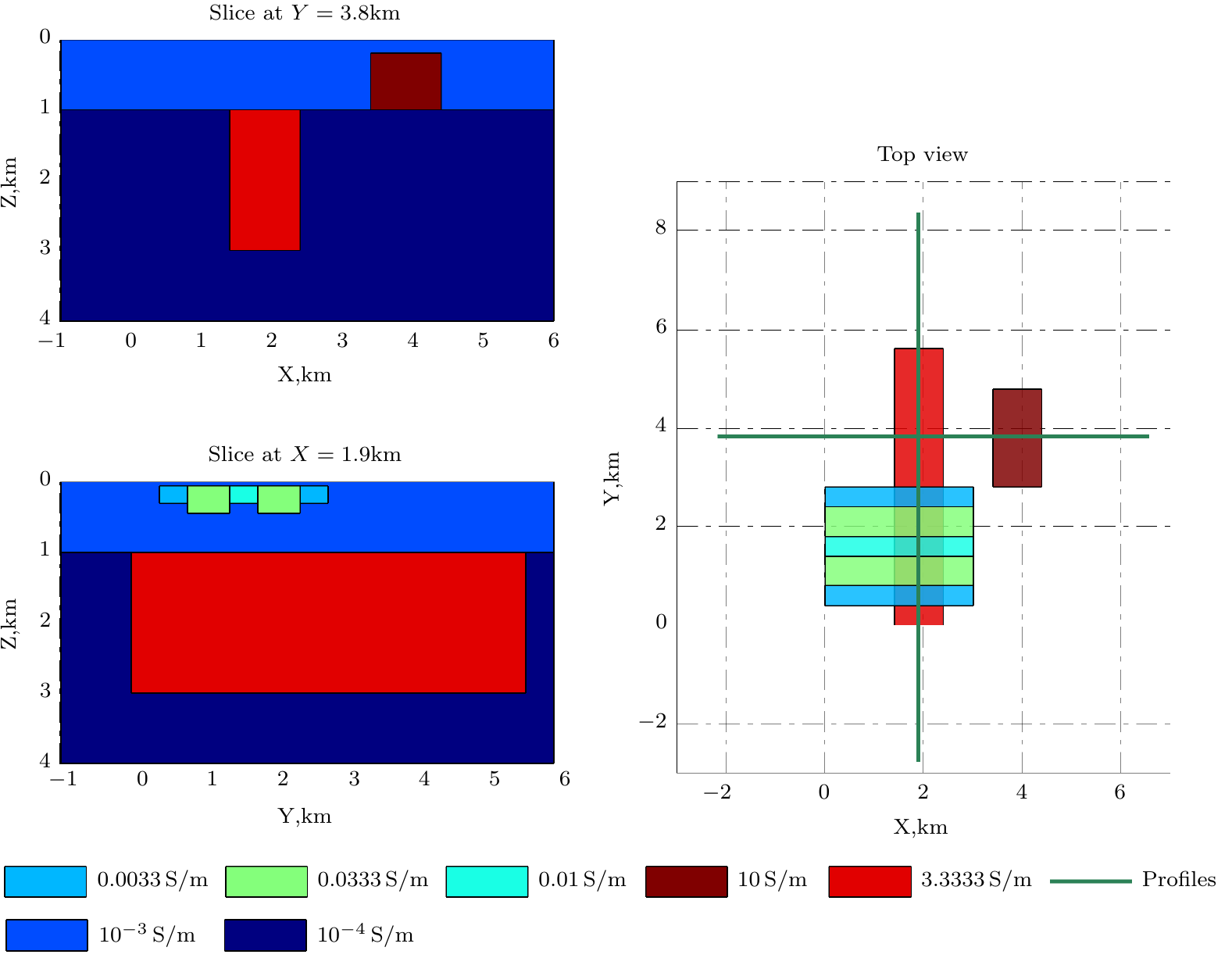}
	\caption{Model COMMEMI3D-3}
	\label{commemi3d3_model}
\end{figure}

The COMMEMI3D-3 model consists of seven rectangular blocks placed in a layered media and oriented along coordinate axes. Their conductivities $\sigma$ (in S/m) and positions (coordinates of the opposite corners in km) are listed in  Table~\ref{commemi_table} and depicted in~Fig.~\ref{commemi3d3_model}.
The layered background  of the model consists of the upper halfspace~${z<0}$~(air) with conductivity of $0$ S/m, two layers with the  conductivity of $10^{-3}$,~$10^{-4}$~S/m, and a lower halfspace with  conductvity~of~$0.1$~S/m. The thickness of the first and the second layers is 1 and 6.5\,km respectively. One can see that maximum conductivity contrast  is $10^4$ in the first layer and $3.3 \cdot 10^4$ in the second one. 

\begin{figure}
	\includegraphics[width=\textwidth , keepaspectratio]{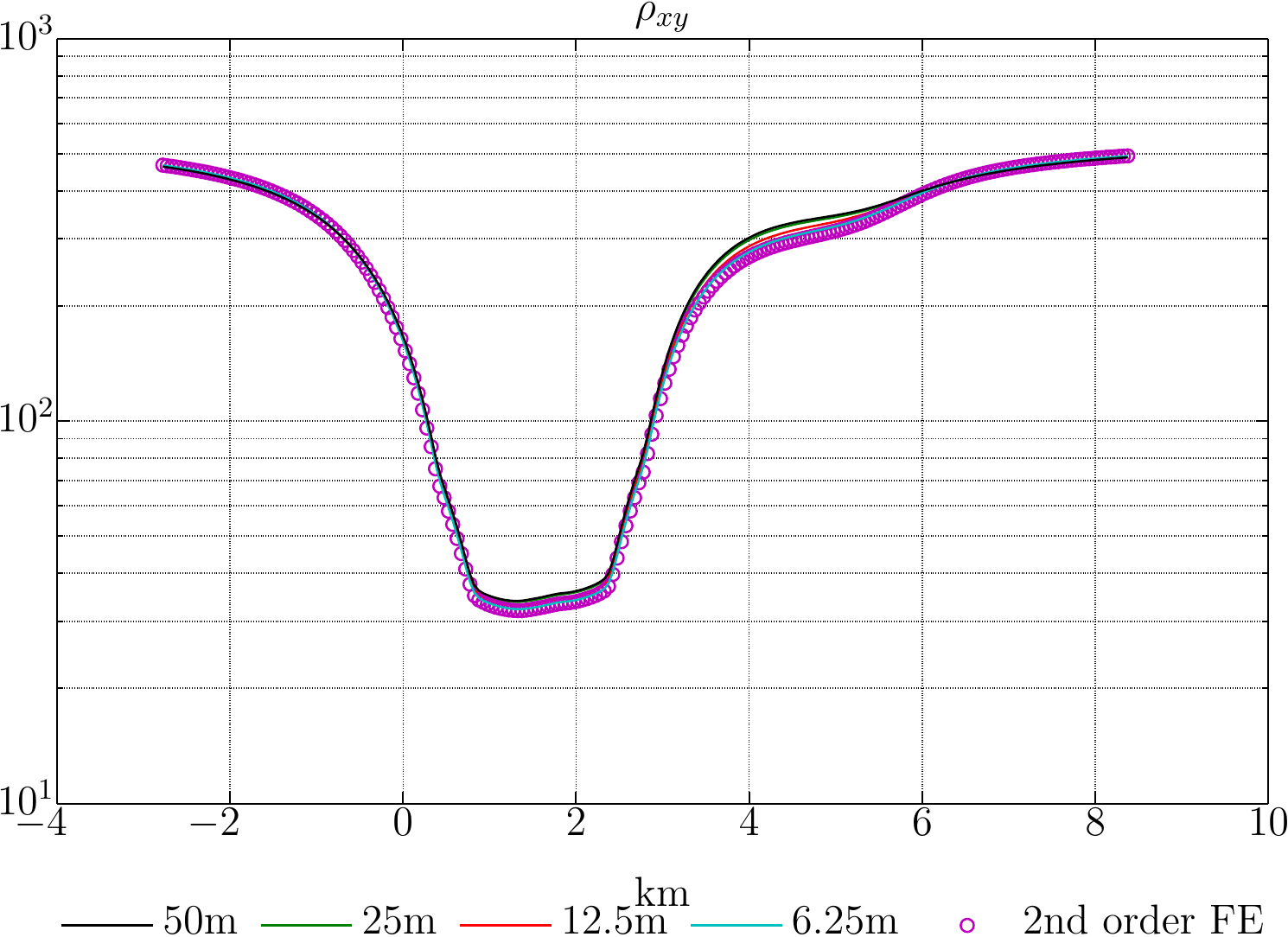}
	 \caption {Apparent resistivity $\rho_{xy}$ at period 1\,s along profile $x=1.9$\,km  for different  subdomain sizes}
	\label{results_rho_xy_x1900}
\end{figure}

\begin{figure}
	\includegraphics[width=\textwidth , keepaspectratio]{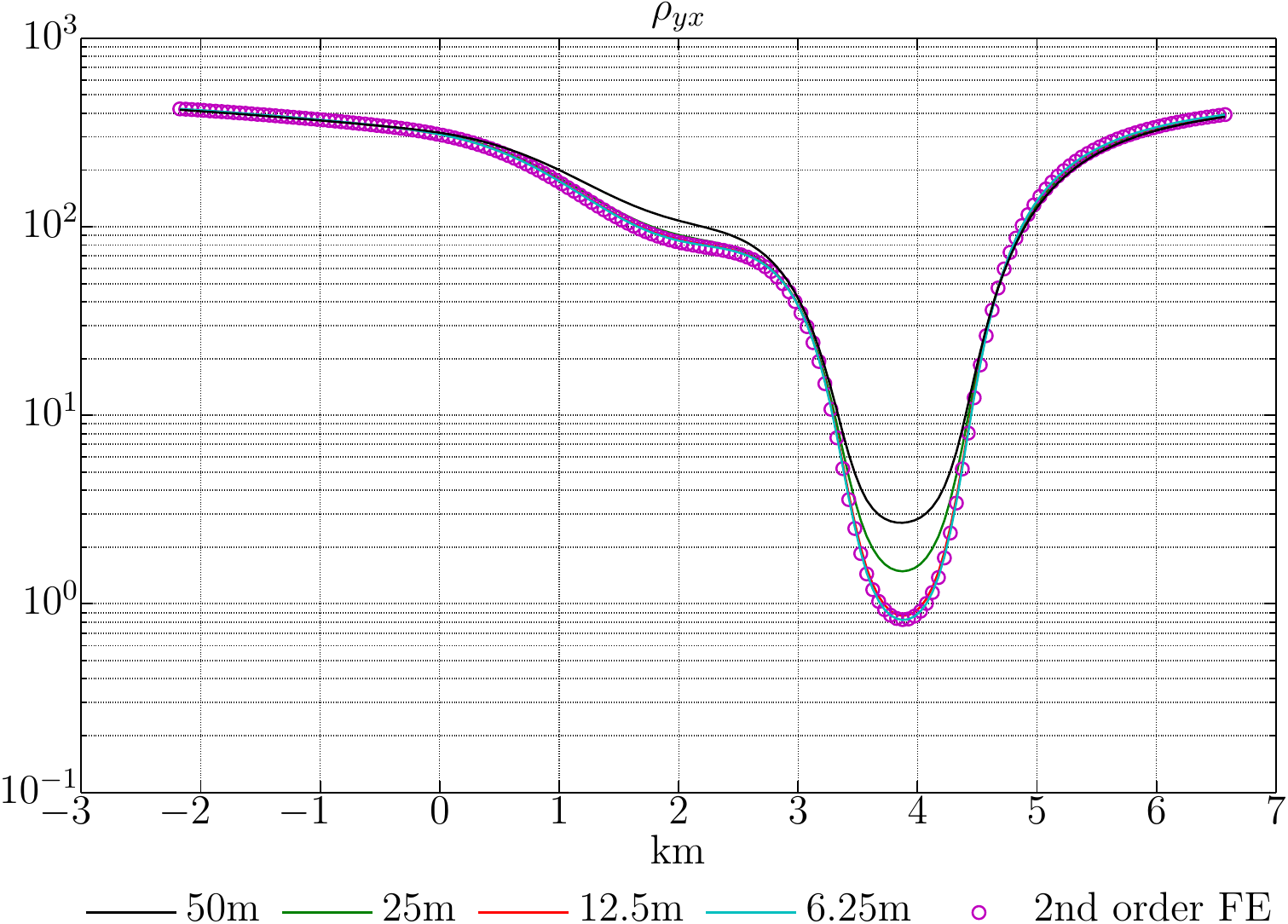}
	 \caption {Apparent resistivity $\rho_{yx}$ at period 1\,s along profile $y=3.83$\,km  for different  subdomain sizes}
	\label{results_rho_yx_y3830}
\end{figure}

\begin{table}
\begin{center}
	\begin{tabular}{|c|c|c|c|c|c|c|c|}\hline
		& $x_1$&$y_1$  & $z_1$   & $x_2$ & $y_2$   & $z_2$  & $\sigma$ \\ 
		\hline
		1& 0   &2.4 & 0.05 & 3   & 2.8  & 0.3  & 0.0033 \\ \hline
	        2& 0   &1.8 & 0.05 & 3   & 2.4  & 0.45 & 0.033  \\ \hline
	        3& 0   &1.4 & 0.05 & 3   & 1.8  & 0.30 & 0.1    \\ \hline
	        4& 0   &0.8 & 0.05 & 3   & 1.4  & 0.45 & 0.033  \\ \hline
	        5& 0   &0.4 & 0.05 & 3   & 0.8  & 0.30 & 0.0033 \\ \hline
	        6& 3.4 &2.8 & 0.2  & 4.4 & 4.8  & 1    & 10     \\ \hline
	        7& 1.4 &0   & 1    & 2.4 & 5.6  & 3    & 3.3333 \\ \hline
	\end{tabular}
\end{center}
\caption{The coordinates of the opposite corners $(x_1,y_1,z_1)$, $(x_2,y_2,z_2)$ in km and conductivities~$\sigma$ (S/m) of COMMEMI3D-3 blocks}
\label{commemi_table}
\end{table}
The modeling of magnetotelluric sounding  was performed for various discretizations (cubic subdomains of different sizes) and periods $T=2\pi/\omega$,  and was compared with the results from modern FE solver by \cite{Grayver2015}.
Figures \ref{results_rho_xy_x1900} and~\ref{results_rho_yx_y3830} represent, correspondingly, the apparent resistivities $\rho_{xy}$ at profile $x=1.9$\,km and $\rho_{yx}$ at profile $y=3.83$\, km for the period 1\,s. One can see that the agreement with FE (magenta circles) is good  even for rather coarse anomaly discretization (black curve). The exception is the area  $[3.5, 4.5]$\,km on the profile $y=3.83$\,km (Fig. \ref{results_rho_yx_y3830}) above the high-conductivity block. This is amended by taking finer discretization (azure curve).

\begin{figure}
	\includegraphics[width=\textwidth , keepaspectratio]{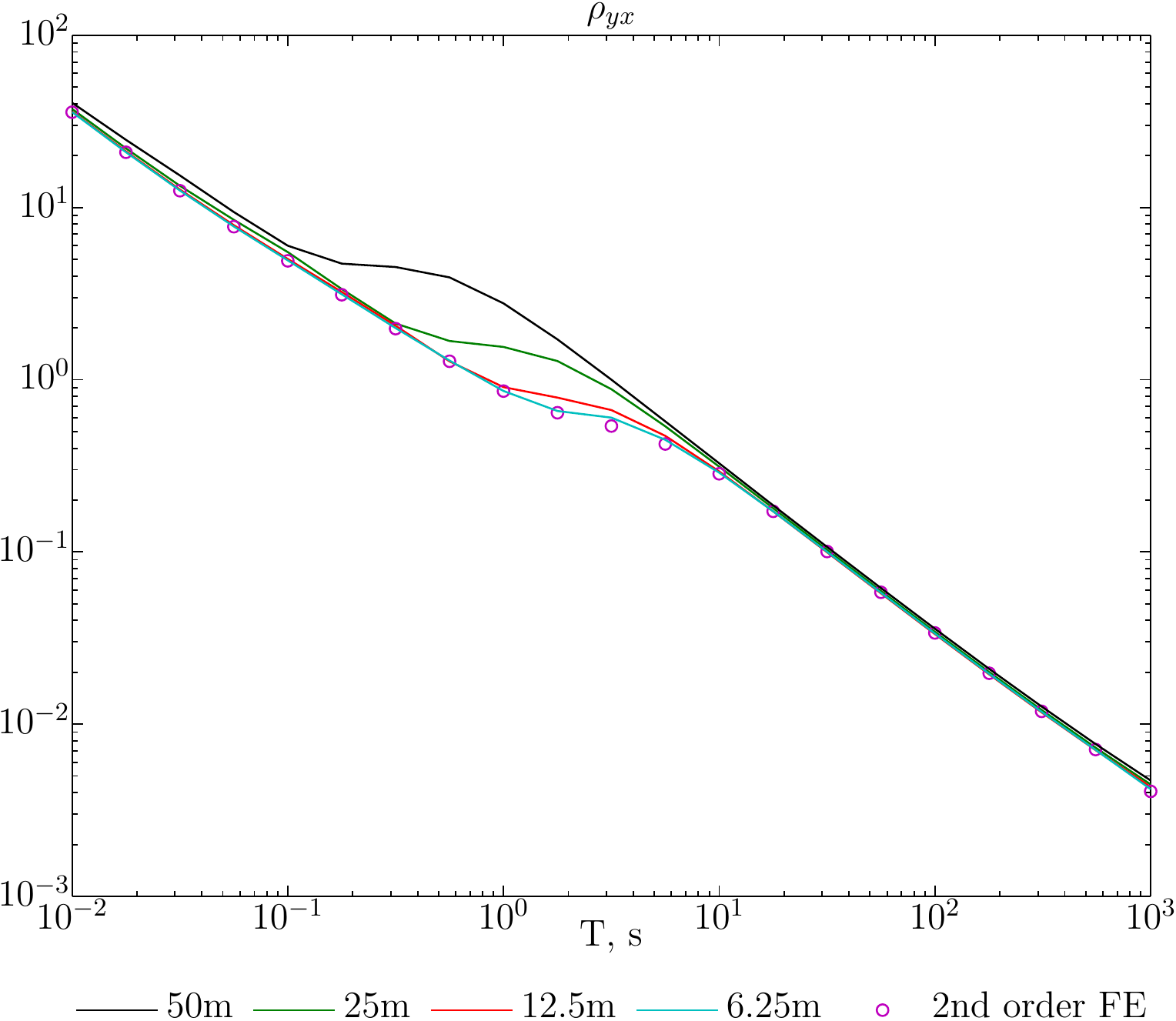}
	 \caption {Apparent resistivity $\rho_{yx}$  at site $x=3.975$\,km, $y=3.83$\,km   depends on period}
	 \label{results_rho_yx_freq}
\end{figure}

\begin{figure}
	\includegraphics[width=\textwidth , keepaspectratio]{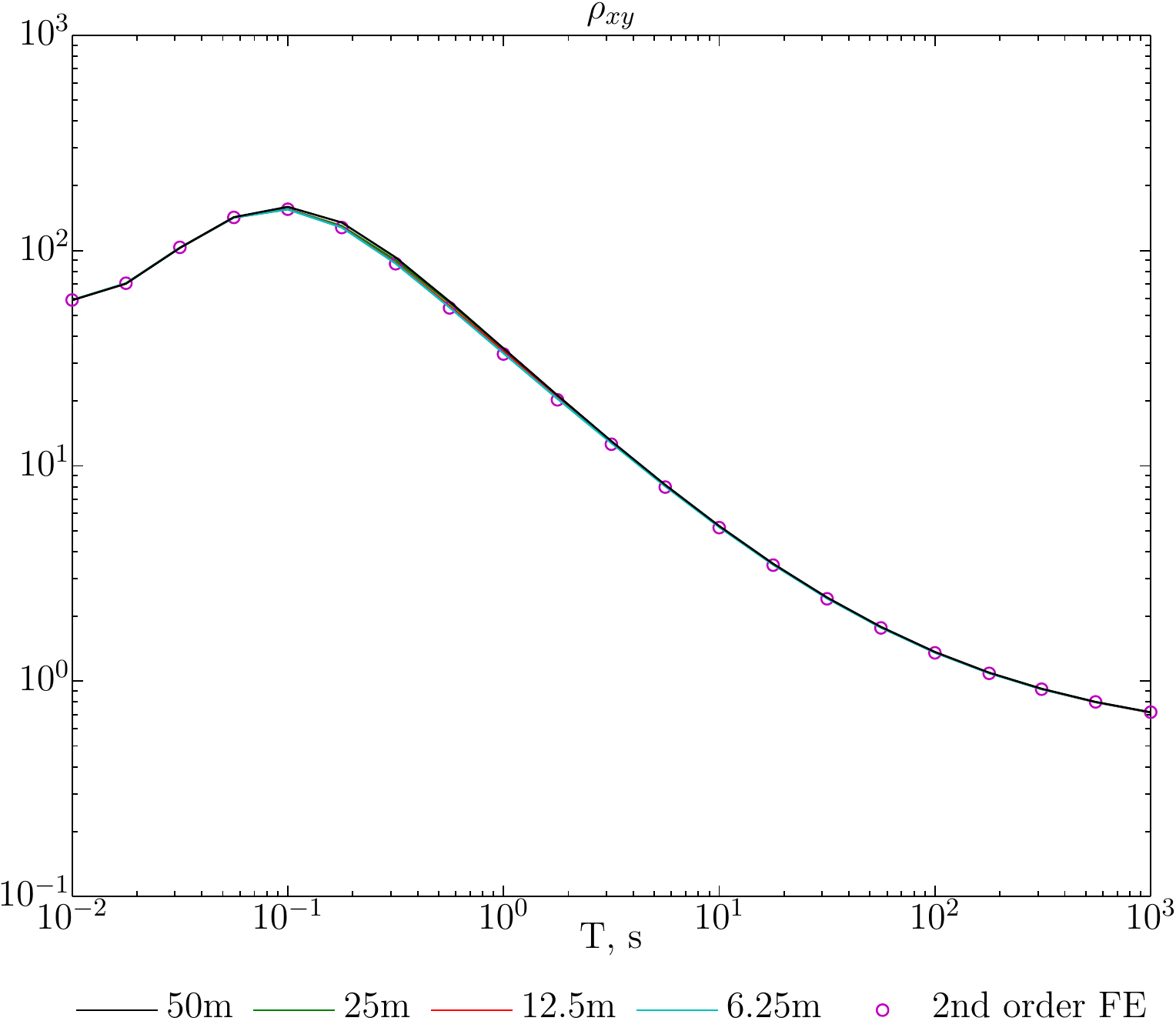}
	 \caption {Apparent resistivity $\rho_{xy}$  at site $x=1.9$\,km, $y=1.7$\,km   depends on period}
	 \label{results_rho_xy_freq}
\end{figure}

Figure \ref{results_rho_yx_freq} shows the apparent conductivity $\rho_{yx}$ at site $(3.975,3.83)$\, km, that is above the high-conductivity block, at different periods.
One can see that the finer discretization is needed only for the periods $\left[10^{-1}, 10^1\right]$\,s.
The reason behind this effect is the drastic change in the electric field inside of the compact high-conductivity block that does not allow to use the piecewise approximation on coarse discretization. The solver by \cite{Grayver2015} uses the second order polynomials which are very effective in such situations.
At the same time, the coarse discretization can be efficiently used for smaller and larger periods. It is worth emphasizing that this concerns only the area above the high-conductivity block, while at the point $(1.9,1.7)$\,km, Fig.~\ref{results_rho_xy_freq} shows good correspondence for all periods.

\begin{figure}
	\includegraphics[width=\textwidth , keepaspectratio]{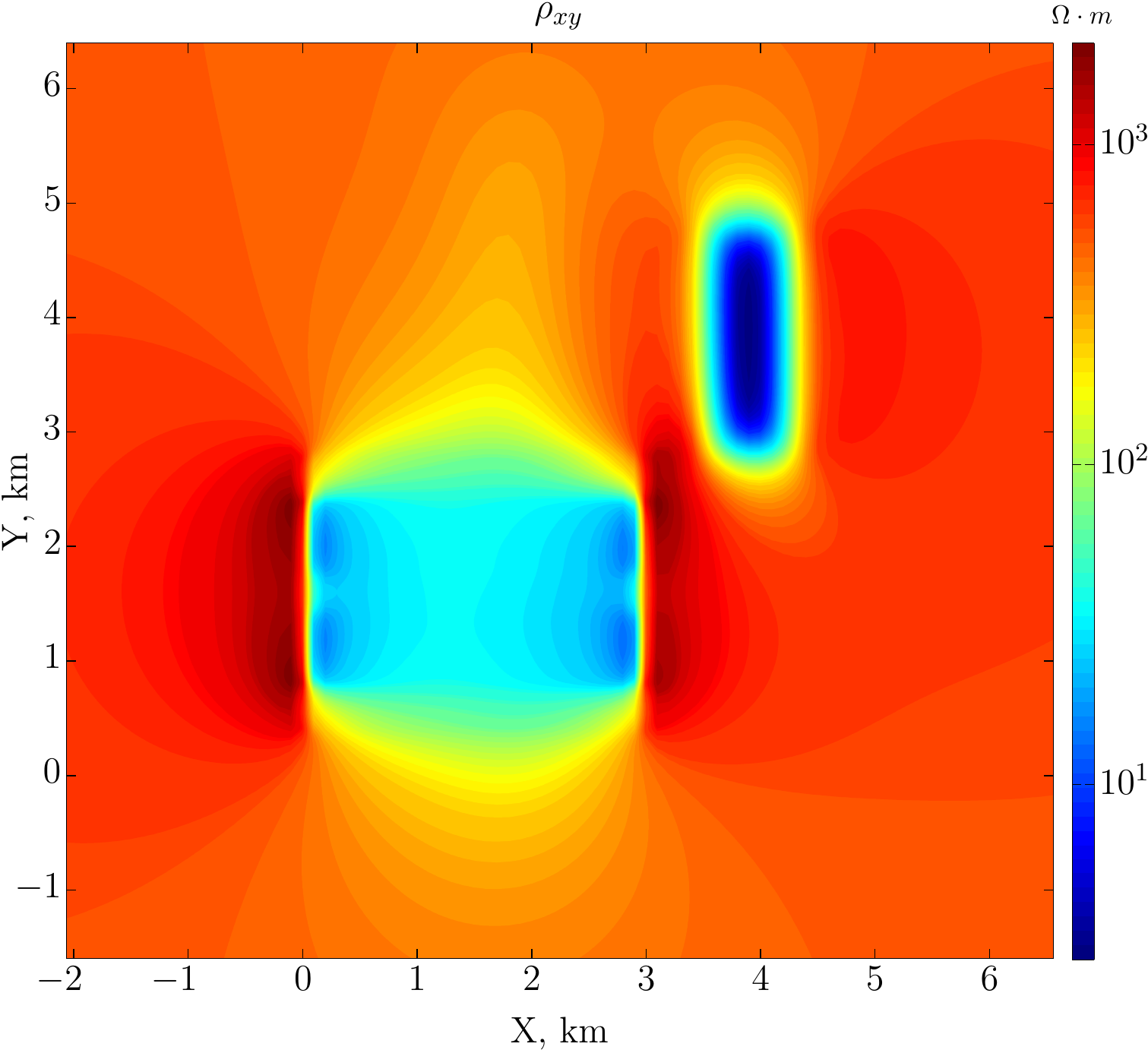}
	 \caption {Apparent resistivity $\rho_{xy}$ at  1\,s}
	\label{results_rho_xy}
\end{figure}

\begin{figure}
	\includegraphics[width=\textwidth , keepaspectratio]{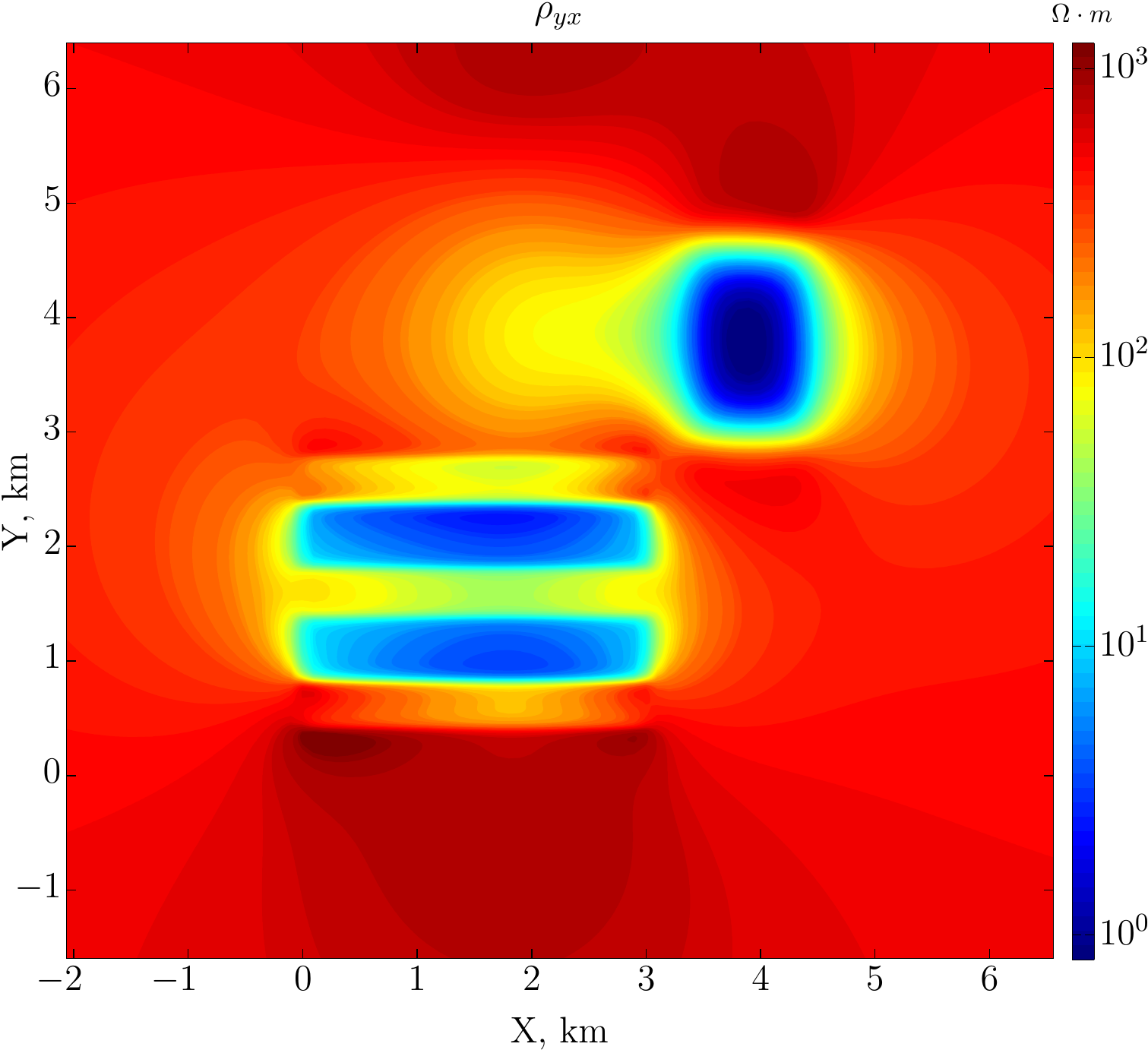}
	 \caption {Apparent resistivity $\rho_{yx}$ at  1\,s}
	\label{results_rho_yx}
\end{figure}

\begin{figure}
	\includegraphics[width=\textwidth , keepaspectratio]{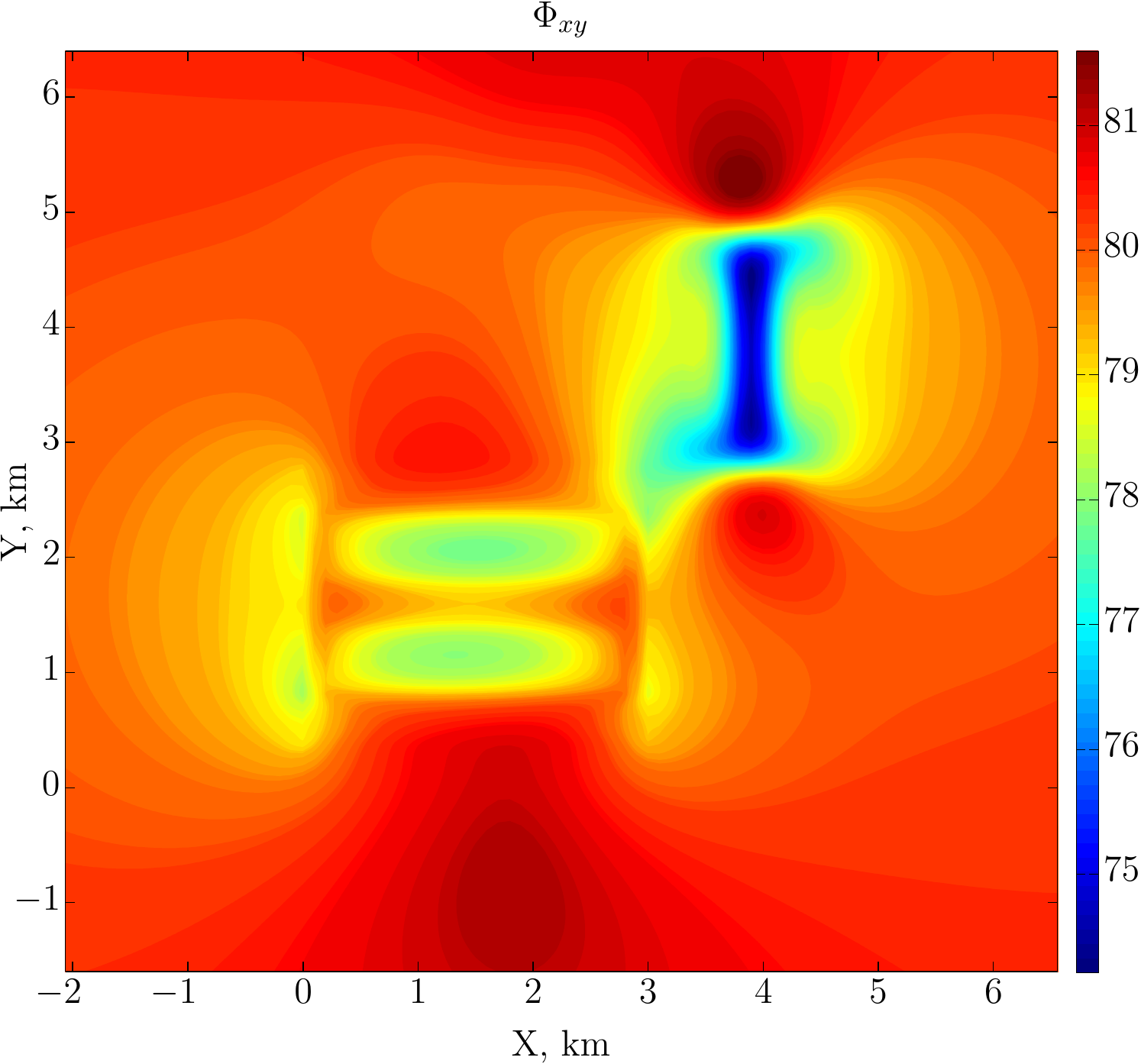}
	 \caption {Impedance phase  $\Phi_{xy}$  at 1\,s }
	\label{results_phi_xy}
\end{figure}

\begin{figure}
	\includegraphics[width=\textwidth , keepaspectratio]{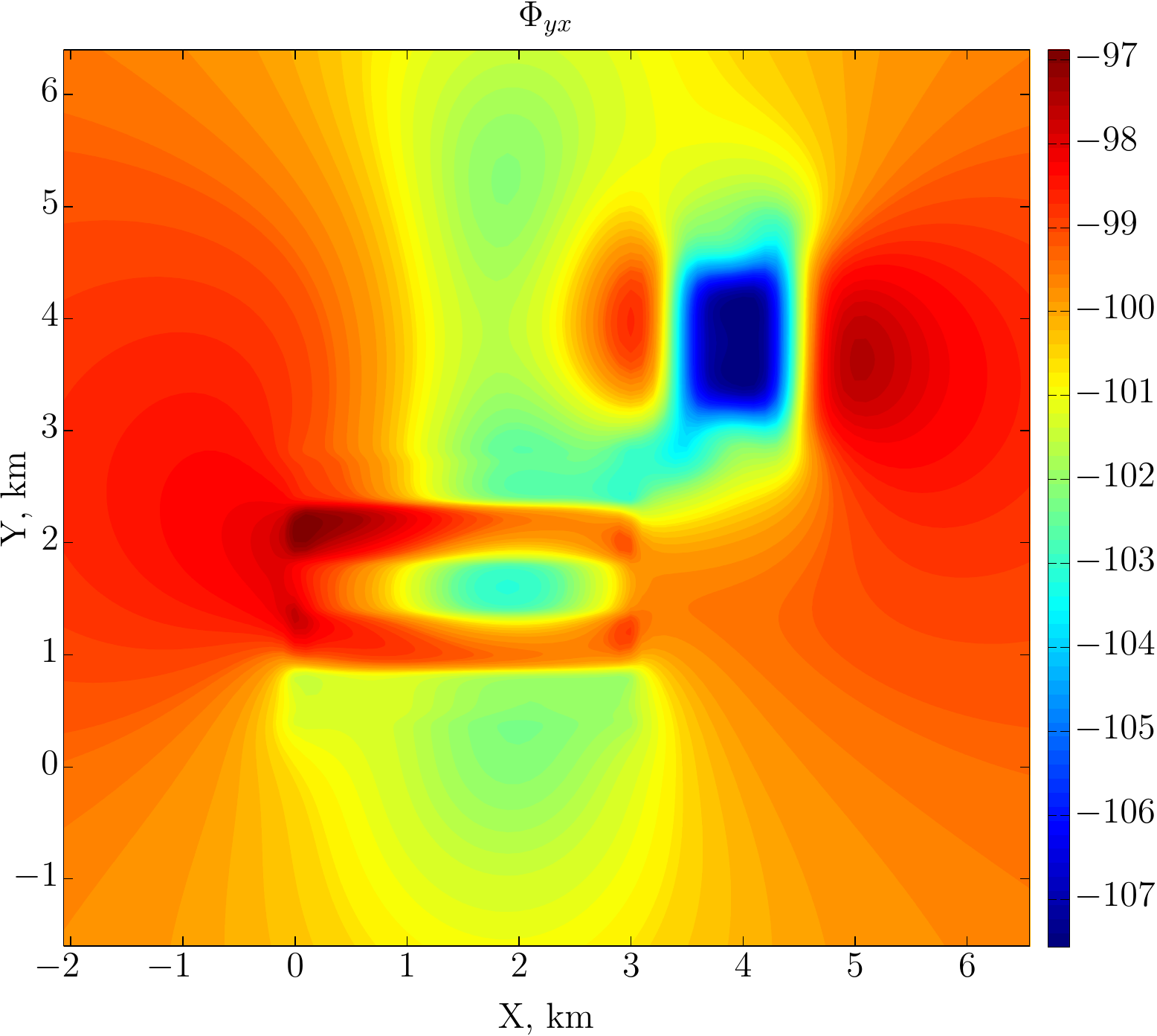}
	 \caption {Impedance phase  $\Phi_{yx}$  at 1\,s}
	\label{results_phi_yx}
\end{figure}

Figures \ref{results_rho_xy} to \ref{results_phi_yx} demonstrate area distribution of the apparent conductivity and impedance phases for 1\,s period. One can see that the variation in apparent conductivity is of the four order of magnitude with very drastic transition. The phase change of the impedance  is 10 degrees, and the transition is again drastic. It is worth reminding, that one of the peculiarities of IE method is quite weak dependence of computational costs on the number of sites where the field is computed. This allowed to obtain the maps with such drastic transitions without the additional computational costs.




\section{Conclusion}

 The presented solver named  ``Gnu Integral Equation Modeling in ElectroMagnetic Geophysics''  (GIEM2G)  shows impressive performance in terms of both memory requirements and  accuracy.
  The memory requirements are 8 times lower compared to other volumetric IE solvers \cite{AvdeevEtAl1997}, \cite{Hursan02}.
 It is achieved for any layered background and non uniform discretization in vertical direction. In this way the average-scale modeling (up to $3\cdot 10^6$ subdomains) can be efficiently done using  laptops. The parallelization scheme allows to use HPC with hundreds and thousands of nodes  for large-scale modeling (up to $10^9$ subdomains).

 The computational efficiency of the  method is demonstrated on high-conductive contrast ($3.3 \cdot 10^4$) model COMMEMI3D-3. To the best of the authors knowledge, it is the first time that 
 such  high-contrast complex model provides comparable results for such different methods as 
 FE and IE. 
 In addition to the efficient usage of HPC the proposed IE method relies on the new 
technique to calculate the matrix coefficients.
It is based on the analytical integration in vertical direction and completely new scheme to compute the integrals in horizontal direction.
It is worth mentioning that the proposed scheme of analytical integration is robust in terms of machine precision and needs only $O(N_z)$ computations of complex exponents.

GIEM2G is implemented  as hybrid MPI+OpenMP software on modern Fortran language. It is an open source software  distributed under the  GPLv2 license and can be simple cloned from GitLab by {\textit{git -clone git@gitlab.com:m.kruglyakov/GIEM2G.git}}. It is also used as an optional computational engine in forward solver extrEMe \cite{Kruglyakov2016} and inverse solver extrEMe-I \cite{KruglyakovAGU2016}.
 

\addcontentsline{toc}{section}{Acknowledgements}
\begin{acknowledgements}
The research of the first author was supported  by the Russian Foundation for Basic Research (grant no. 13-05-12018-OFI\_M). As a visiting fellow in ETH Zurich he was also partially supported by the Swiss National Science Foundation (grant no. IZK0Z2\_163494) and ETH Zurich.
Authors acknowledge the teams of  HPC CMC Lomonosov MSU  for the access to ``Bluegene/P'' HPC,  the  Lomonosov MSU Research Computing Center for the access to HPC ``Lomonosov'' \cite{HPCLomonosov} and the Swiss National
 Supercomputing Center (CSCS) grant (project ID s577). Authors also would like to thank Alexander Grayver, ETH Zurich, for providing data for comparison and  Alexey Kuvshinov, ETH Zurich, for suggestions and helpful discussions.
\end{acknowledgements}


\bibliographystyle{spbasic}   
\bibliography{giem2g_bib}
\addcontentsline{toc}{section}{References}

\begin{appendices}
\renewcommand{\thesection}{\Alph{section}:}
\makeatletter
\newcommand{\applabel}[1]{\def\@currentlabel{\Alph{section}}\label{#1}}
\makeatother

\section{Green's Tensor}\applabel{AppendixGreensTensor}

Following the notations from \cite{dmitriev2002}, the electrical  $\GE$ and magnetic~$\GH$ tensors of layered media can be written as
\begin{equation}
\label{geh}
	\begin{aligned}
		\GE&={\widehat G+\grad \left(\frac{\mu_0}{k^2}\DIV\frac{\widehat G}{\mu_0}\right)},\\
		\GH&=\frac{1}{i\omega\mu_0}\rot \widehat G,\\
		k^2&=i\omega\mu_0\tsig_b,
	\end{aligned}
\end{equation}
where
\begin{equation}
\label{green_tensor}
	\widehat G(M,M_0)=\begin{pmatrix}
		G_1(M,M_0) & 0 & 0\\
		0 & G_1(M,M_0) & 0\\
		\prt{g(M,M_0)}{x} & \prt{g(M,M_0)}{y} & G_2(M,M_0)
	\end{pmatrix}
\end{equation}
and
	\begin{align}
	\label{G1G2g}
		&G_1(M,M_0)=\frac{i\omega\mu_0}{4\pi}\int\limits_0^\infty J_0(\lambda\rho)U_1(\lambda,z,z_0)\lambda d\lambda,\notag\\
		&G_2(M,M_0)=\frac{i\omega\mu_0}{4\pi}\int\limits_0^\infty J_0(\lambda\rho)U_\sigma(\lambda,z,z_0)\lambda d\lambda,\notag\\
		&g(M,M_0)=-\frac{i\omega\mu_0}{4\pi}\int\limits_0^\infty J_0(\lambda\rho)\left(\prt{ }{z_0}U_\sigma(\lambda,z,z_0)+\prt{ }{z} U_1(\lambda,z,z_0)\right)\frac{d\lambda}{\lambda},\notag\\
		&M=M(x,y,z)~ ~ ~M_0=M_0(x_0,y_0,z_0)~ ~ ~\rho=\sqrt{(x-x_0)^2+(y-y_0)^2}.
	\end{align}
Here $J_0$ is a zero-order Bessel function of the first kind, and functions $U_\gamma(\lambda,z,z_0)$, $\gamma=1,\sigma$ are the fundamental functions of the layered media (Appendix~\ref{AppendixVerticalIntegration}).

From \eqref{geh} and \eqref{green_tensor} one gets 
\begin{equation}
 \label{Ge2}
  \GE =\left\{\begin{matrix}
	  G_1+\frac{1}{k^2}\prtt{ }{x}\left(G_1+\prt{g}{z}\right)& \frac{1}{k^2}\prtq{ }{x}{y}\left(G_1+\prt{g}{z}\right)&\frac{1}{k^2}\prtq{G_2}{x}{z}\\
\frac{1}{k^2}\prtq{ }{x}{y}\left(G_1+\prt{g}{z}\right)& G_1+\frac{1}{k^2}\prtt{ }{y}\left(G_1+\prt{g}{z}\right)& \frac{1}{k^2}\prtq{G_2}{y}{z}\\
\prt{g}{x}+\frac{1}{k^2}\prtq{ }{x}{z}\left(G_1+ \prt{g}{z}\right)&\prt{g}{y}+\frac{1}{k^2}\prtq{ }{y}{z}\left(G_1+\prt{g}{z} \right)& G_2+\frac{1}{k^2}\prtt{G_2}{z}
                   \end{matrix}\right\}.
\end{equation}
Note, that $G^E_{xz}(M,M_0)=-G^E_{zx}(M_0,M)$, $G^E_{yz}(M,M_0)=-G^E_{zy}(M_0,M)$ according to Lorentz reciprocity.

Let $\AD_n=S_n \times [z_n^1,z_n^2]$, $\AD_m=S_m \times [z_m^1,z_m^2]$, where $S_n, S_m$ are horizontal rectangular domains, and let
$\tsig_b=\tsig_b^n$ inside $[z_n^1, z_n^2]$,  $\tsig_b=\tsig_b^m$ inside $[z_m^1, z_m^2]$. That is the subdomains do not intersect the boundaries of the layers. Taking into account \eqref{slae_matricies}, \eqref{geh} and \eqref{Ge2} one can see that $\widehat B_n^m$ is expressed in terms of  double volumetric  integrals  with weak integrable singularity, so the order of integration can be changed. Then using  \eqref{G1G2g} and \eqref{Ge2}  one obtains
\begin{equation} 
	\widehat B_n^m=\frac{1}{4\pi}\left\{\begin{aligned}
		&I_1^{n,m}+I_{xx}^{n,m} && I_{xy}^{n,m}          &&& I_x^{n,m}\\
		& I_{xy}^{n,m}          &&I_1^{n,m}+I_{yy}^{n,m} &&& I_y^{n,m}\\
		&- I_x^{m,n}            &&- I_y^{m,n}            &&&I_2^{n,m}+I_{zz}^{n,m}\\
	\end{aligned}\right\},
\end{equation}
where
\begin{align}
    I_1^{n,m}&=\int\limits_{S_n}\int\limits_{S_m}\left(\int\limits_0^{\infty}J_0(\lambda\rho)V_1^{n,m}(\lambda)\lambda d\lambda\right) dx_0dy_0dxdy,\notag\\
    I_2^{n,m}&=\int\limits_{S_n}\int\limits_{S_m}\left(\int\limits_0^{\infty}J_0(\lambda\rho)V_2^{n,m}(\lambda)\lambda d\lambda\right) dx_0dy_0dxdy,\notag\\
    I_{\alpha\beta}^{n,m}&=\int\limits_{S_n}\prtq{ }{\alpha}{\beta}\left\{\int\limits_{S_m}\left(\int\limits_0^{\infty}J_0(\lambda\rho)\left[V_1^{n,m}(\lambda)+V_3^{n,m}(\lambda)\right]\frac{d\lambda}{\lambda}\right) dx_0dy_0\right\}dxdy,\notag\\
    I_{\alpha}^{n,m}&=\int\limits_{S_n}\prt{}{\alpha}\left\{\int\limits_{S_m}\left(\int\limits_0^{\infty}J_0(\lambda\rho)V_4^{n,m}(\lambda)\lambda d\lambda\right) dx_0dy_0\right\}dxdy,\notag\\
    I_{zz}^{n,m}&=\int\limits_{S_n}\int\limits_{S_m}\left(\int\limits_0^{\infty}J_0(\lambda\rho)V_5^{n,m}(\lambda)\lambda d\lambda\right) dx_0dy_0dxdy.\label{Hor_Int}
\end{align}
Here $\alpha=x,y$, $\beta=x,y$, and 
\begin{flalign}
    &V_1^{n,m}(\lambda)=i\omega\mu_0\int\limits_{z_n^1}^{z_n^2}\int\limits_{z_m^1}^{z_m^2}U_1(\lambda,z,z_*)dz_*dz,&&\notag\\
    &V_2^{n,m}(\lambda)=i\omega\mu_0\int\limits_{z_n^1}^{z_n^2}\int\limits_{z_m^1}^{z_m^2}U_\tsig(\lambda,z,z_*)dz_*dz,&&\notag\\
    &V_3^{n,m}(\lambda)=-\frac{1}{\tsig_b^n}\int\limits_{z_n^1}^{z_n^2}\prt{ }{z} \left(\int\limits_{z_m^1}^{z_m^2}\left[\prt{ }{z_*}U_\tsig(\lambda,z,z_*)\right]dz_*\right) dz,\notag\\    
    &V_4^{n,m}(\lambda)=\frac{1}{\tsig_b^n}\int\limits_{z_n^1}^{z_n^2}\prt{ }{z}\left(\int\limits_{z_m^1}^{z_m^2}U_\tsig(\lambda,z,z_*)dz_*\right)dz,&&\notag\\
    &V_5^{n,m}(\lambda)=\frac{1}{\tsig_b^n}\int\limits_{z_n^1}^{z_n^2}\prtt{ }{z}\left(\int\limits_{z_m^1}^{z_m^2}U_\tsig(\lambda,z,z_*)dz_*\right)dz.&&\label{Ver_Int}
\end{flalign}
Therefore, to obtain the  coefficients of $\hat B_n^m$ one needs computational methods to find ``horizontal'' integrals \eqref{Hor_Int} and ``vertical'' integrals \eqref{Ver_Int}. These methods are presented in the next sections.


\section{Vertical Integration}\applabel{AppendixVerticalIntegration}
The integrals in \eqref{Ver_Int} are expressed in terms of the so-called fundamental function of the layered media, \cite{dmitriev2002}. 
Consider the media with $N_{lay}-1$  homogenous layers with complex conductivities $\tsig_n$, $n=1\dots N_{lay}-1$, the upper halfspace (air, the zeroth layer) with complex conductivity $\tsig_0$ and the lower halfspace (the $N_{lay}$-th layer) with conductivity $\tsig_{N_{lay}}$.
Note, that in EM sounding problems typically $\Re{\tsig_0} \le 10^{-9}$.

The  function $U_\gamma(z,z_*,\lambda)$ is defined as a unique  solution of the problem
\begin{equation}
	\label{lmff}
	\left\{
	\begin{aligned}
		&\prt{^2 }{z^2}U_\gamma(z,z_*,\lambda) -\eta_0^2U_\gamma(z,z_*,\lambda)=0,z <d_1,  z \not =z_*,\\
		&\prt{^2 }{z^2}U_\gamma(z,z_*,\lambda) -\eta_n^2U_\gamma(z,z_*,\lambda)=0, d_n<z<d_{n+1}, z \not =z_*,\\
		&\prt{^2 }{z^2}U_\gamma(z,z_*,\lambda) -\eta_{N_{lay}}^2U_\gamma(z,z_*,\lambda)=0, z >d_{N_{lay}}, z \not =z_*,\\
		&\left[U_\gamma(z,z_*,\lambda)\right]_{z=d_n}=0,\\
		&\left[\frac{1}{\gamma}\prt{ }{z}U_\gamma(z,z_*,\lambda)\right]_{z=d_n}=0,\\
		&\left[U_\gamma(z,z_*,\lambda)\right]_{z=z_*}=0,\\
		&\left[\prt{ }{z}U_\gamma(z,z_*,\lambda)\right]_{z=z_*}=-2,\\
		&\left|U_\gamma(z,z_*,\lambda)\right| \to  0 \text{~as~}z \to \pm \infty,\\
		&\eta_m^2=\lambda^2-k^2_m,\ ~ ~ ~k^2_m=i\omega\mu_0\tsig_m,\\
		&m=0\dots N_{lay},\quad n=1\dots N_{lay}-1,\quad 0< \lambda < \infty.
	\end{aligned}
	\right.
\end{equation}


The following procedure is performed to obtain an  explicit expression for $U_\gamma$ that allows analytical integration.
 Let $l_0=0$, $l_{N_{lay}}=0$, $l_n=d_{n+1}-d_n$, $n=1\dots N_{lay}-1$. Define $p_{m}^\gamma$, $q_{m}^\gamma$, ${m=0\dots N_{lay}}$ by the recurrent expressions
\begin{equation}
	\label{pq_req}
	\begin{aligned}
		&p_0^\gamma=0; &&q_{N_{lay}}^\gamma=0;\\ 
		&p_{1}^\gamma=\frac{1-\alpha_0^\gamma\frac{\eta_0}{\eta_{1}}}{1+\alpha_0^\gamma \frac{\eta_0}{\eta_{1}}};\quad
		&&q_{N_{lay}-1}^\gamma=\frac{1-\beta_{N_{lay}}^\gamma\frac{\eta_{N_{lay}}}{\eta_{N_{lay}-1}}}{1+\beta_{N_{lay}}^\gamma\frac{\eta_{N_{lay}}}{\eta_{N_{lay}-1}}};\\
	&p_{m+1}^\gamma=\frac{1+\alpha_m^\gamma\frac{\eta_m}{\eta_{m+1}}\frac{p_{m}^\gamma e^{-2\eta_ml_m}-1}{p_{m}^\gamma e^{-2\eta_ml_m}+1}}{1-\alpha_m^\gamma \frac{\eta_m}{\eta_{m+1}}\frac{p_{m}^\gamma e^{-2\eta_ml_m}-1}{p_{m}^\gamma e^{-2\eta_ml_m}+1}}, \quad 
	&&q_{m-1}^\gamma=\frac{1+\beta_m^\gamma\frac{\eta_m}{\eta_{m-1}}\frac{q_{m}^\gamma e^{-2\eta_ml_m}-1}{q_{m}^\gamma e^{-2\eta_ml_m}+1}}{1-\beta_m^\gamma \frac{\eta_m}{\eta_{m-1}}\frac{q_{m}^\gamma e^{-2\eta_ml_m}-1}{q_{m}^\gamma e^{-2\eta_ml_m}+1}},\\
	&\hspace{1.8cm} m\not =N_{lay}; \qquad &&\hspace{1.8cm} m\not =0;\\
	&\alpha_m^\gamma=\left\{\begin{aligned}
		1,\gamma=1;\\
		\frac{\tsig_{m+1}}{\tsig_{m}},\gamma=\tsig;
	\end{aligned}
	\right.
	&&\beta_m^\gamma=\left\{\begin{aligned}
		1,\gamma=1;\\
		\frac{\tsig_{m-1}}{\tsig_{m}},\gamma=\tsig.
	\end{aligned}
	\right.
\end{aligned}
 \end{equation}
Let $d_0=d_1$, $d_{N_{lay+1}}=d_{N_{lay}}$ and let points $z_r$, $z_s$ belong to $r$ and $s$ layers respectively, $0 \le r, s \le N_{lay}$. Then using \eqref{pq_req} one gets  
\begin{equation}\label{U_explicit}
 U_\gamma(z_r,z_s,\lambda)=
 \begin{cases}
              &A_{r,s}^\gamma\left(p_r^\gamma e^{2\eta_rd_r}e^{-\eta_rz_r}+e^{\eta_r z_r}\right)\left(e^{-\eta_sz_s}+q_s^\gamma e^{-2\eta_sd_{s+1}}e^{\eta_s z_s}\right)\\
	&\hspace{6.cm}\text{for}\quad z_r\le z_s;\\
		&A_{s,r}^\gamma\left(p_s^\gamma e^{2\eta_sd_s}e^{-\eta_sz_s}+e^{\eta_s z_s}\right)\left(e^{-\eta_rz_r}+q_r^\gamma e^{-2\eta_rd_{r+1}}e^{\eta_r z_r}\right)\\
		&\hspace{6.cm}\text{for}\quad z_r>z_s,              
 \end{cases}
\end{equation}
where
\begin{equation}
	\label{A_rs}
	\begin{aligned}
		A_{r,s}^\gamma&=Q_r^\gamma\times Q_{r+1}^\gamma\times \cdots \times Q_{s-1}^\gamma A_{s,s}^\gamma,\quad \text{ for\quad $r < s$},\\
		Q_m^\gamma&=\frac{1+p_{m+1}^\gamma } {1+p_m^\gamma e^{-2\eta_m l_m }}e^{(\eta_{m+1}-\eta_m)d_{m+1}},\quad\text{for}\quad m=1\dots N_{lay}-1,\\
		A_{n,n}^\gamma&=\frac{1}{\eta_n\left(1-p_n^\gamma q_n^\gamma e^{-2\eta_nl_n}\right)},\quad \text{ for\quad $r=s=n$}, n=0\dots N_{lay},\\
		A_{r,s}^1&=A_{s,r}^1,~ ~ ~ A_{r,s}^{\tsig}=\frac{\tsig_r}{\tsig_s}A_{s,r}^{\tsig},\quad \text{ for\quad $ r>s$}.\\
\end{aligned}
\end{equation}
To check \eqref{U_explicit} one can explicitly   substitute \eqref{U_explicit} in \eqref{lmff} taking into account~\eqref{pq_req}~and~\eqref{A_rs}.

In view of \eqref{U_explicit} one can see that integrals in \eqref{Ver_Int} (i.e., the integrals over $U_\gamma(z,z_*,\lambda)$ and its partial derivatives)  can be   integrated analytically with respect to $z$, $z_*$ over any domains that do not  intersect  the layer boundaries. 
However, the rounding errors arising in addition and multiplication of very small or large quantities make the formula~\eqref{U_explicit}   impractical  for $\lambda \gg 1$. Instead the following formula is used 
\begin{equation}\label{UU2}
 U_\gamma(z_r,z_s,\lambda)=
 \begin{cases}
              &A^\gamma_{r,s}\left(p_r^\gamma e^{-(\eta_rz_r+\eta_sz_s-2\eta_rd_r)}+q_s^\gamma e^{-(2\eta_s d_{s+1}-(\eta_rz_r+\eta_sz_s))}\right.+\\
		&\qquad \left.e^{-(\eta_sz_s-\eta_rz_r)}+p_r^\gamma q_s^\gamma e^{-(2(\eta_sd_{s+1}-\eta_rd_r)-(\eta_sz_s-\eta_rz_r))}\right)\\
	&\hspace{6.cm}\text{for}\quad z_r \le z_s;\\
		&A^\gamma_{s,r}\left(p_s^\gamma e^{-(\eta_sz_s+\eta_rz_r-2\eta_sd_s)}+q_r^\gamma e^{-(2\eta_r d_{r+1}-(\eta_sz_s+\eta_rz_r))}\right.+\\
		&\qquad \left.e^{-(\eta_rz_r-\eta_rz_r)}+p_s^\gamma q_r^\gamma e^{-(2(\eta_rd_{r+1}-\eta_sd_s)-(\eta_rz_r-\eta_sz_s))}\right)\\
		&\hspace{6.cm}\text{for}\quad z_r>z_s.              
 \end{cases}
\end{equation}
Formula \eqref{UU2} overcomes the aforementioned problem, since the real parts of all the exponents powers are negative.
 The consequent calculations  provide accurate and robust results for any $0< \lambda <\infty$.

Consider $N_z$ subdomains in the discretization in vertical direction. To obtain the matrix $\widehat B_n^m$ for the system \eqref{bcyl} one needs to compute $O\left(N_z^2\right)$ complex exponents  in~\eqref{UU2}.
An algorithm requiring only $O\left(N_z\right)$ complex exponents calculations is developed to  speed up the integration procedure.

Let $z_0<z_1< \dots< z_{N_z}$. Suppose that the intervals $[z_l,z_{l+1}]$, $l=0\dots N_z-1$ do not intersect the layers' boundaries. 
For $i,j=0\dots N_z-1$, $0 \le \alpha+\beta \le 2$
one needs to calculate  $W_{i,j}^{\alpha,\beta}(\gamma)=\int\limits_{z_i}^{z_{i+1}}\prt{^\alpha }{z^\alpha}\left( \int\limits_{z_j}^{z_{j+1}}\left[\prt{^\beta }{z_*^\beta}\, U_\gamma(z,z_{*},\lambda)\right]dz_{*}\right)dz$. 

Let $r_l$ be an index of the layer containing $[z_l,z_{l+1}]$, $l=0 \dots N_{z}-1$. Then using~\eqref{U_explicit} one obtains for $z_i < z_j$
\begin{equation}
	\label{fact}
	\begin{aligned}
		W_{i,j}^{\alpha,\beta}(\gamma)=&\int\limits_{z_i}^{z_{i+1}}\left(\prt{^\alpha}{z^\alpha} \int\limits_{z_j}^{z_{j+1}}\left[\prt{^\beta }{z_*^\beta} U_\gamma(z,z_*,\lambda)\right]dz_*\right)dz=\\
	&H_\gamma^\alpha(z_i,z_{i+1})\prod\limits_{l={i+1}}^{l={j}}\Theta^\gamma_l \int\limits_{z_j}^{z_{j+1}}\left(\prt{^\beta }{z_*^\beta} U\gamma(z_j,z_*,\lambda)\right)dz_*,\\ 
\end{aligned}
\end{equation}
where
\begin{equation}
	\label{fact_2}
	\begin{aligned}
		H_\gamma^\alpha(z_i,z_{i+1})&=\frac{\int\limits_{z_i}^{z_{i+1}}\left(\prt{^\alpha }{z^\alpha}\left[p_{r_i}^\gamma e^{-\eta_{r_i}\left(z+z_{i+1}-2d_{r_i}\right)}+e^{-\eta_{r_i}\left (z_{i+1}-z\right)}\right]\right)dz}{p_{r_i}^\gamma e^{-2\eta_{r_i}\left(z_{i+1}-d_{r_i}\right)}+1},\\
		\\
		\Theta_l^\gamma&=\varkappa_l^\gamma\frac{ p_{r_{l}}^\gamma e^{-\eta_{r_l}\left((z_l+z_{l+1})-2d_{r_l}\right)}+e^{-\eta_{r_l}( z_{l+1}-z_l)}}{p_{r_{l}}^\gamma e^{2\eta_{r_{l}}(d_{r_{l}}-z_{l+1})}+1},\\
		\varkappa_l^\gamma&=\left\{\begin{aligned}
					1, r_l=r_{l+1}\\
					Q_{r_{l+1}}^\gamma, r_l \not =r_{l+1}\\
					\end{aligned}
		\right\}.\\
	\end{aligned}
\end{equation}

All the exponents   in \eqref{fact_2} vanish as $\lambda \to \infty$, so the corresponding computations do not depend on the round-off errors due to the machine precision. 
The formulas for $z_i>z_j$  are similar.  The integrals  $W_{ii}^{\alpha,\beta}$ are computed analytically using \eqref{UU2}.
Since $\Theta_l^\gamma$ depends only on $l=1\dots N_{z}$ and $\gamma={1,\tsig}$, one only needs to calculate $O\left(N_z\right)$ complex exponents  using factorization \eqref{fact}, \eqref{fact_2}.

 
\section{Horizontal Integration}\applabel{AppendixHorizontalIntegration}

The integrals \eqref{Hor_Int} are the particular case  of the integral
	\begin{equation}
		\label{g_eint}
		\begin{aligned}
			&I_{\alpha,\beta}=\int\limits_{S_n}\prtqq{{\alpha+\beta}}{x^\alpha}{y^\beta}\left\{\int\limits_{S_m}\left[\int\limits_0^\infty J_0(\rho\lambda)f(\lambda)d\lambda\right] dS_m\right\}dS_n,\\
		& \rho=\sqrt{(x-x_0)^2+(y-y_0)^2}, \quad 0 \le \alpha+\beta \le 2,
		\end{aligned}
	\end{equation}
	where $f(\lambda)$ is some easily computed function, $S_n=[x_n,x_n+h_x]\times[y_n,y_n+h_y]$, $S_m=[x_m,x_m+h_x]\times[y_m,y_m+h_y]$ are the rectangular domains with similar sizes.

	The key feature of the proposed method is transformation of integrals \eqref{g_eint} to one-dimensional convolution integral. Taking for simplicity $\alpha=\beta=0$,  one has
	\begin{equation}
	 I_{0,0}=F(R;p,q,\varphi)=\int\limits_0^\infty  K(R\lambda;p,q,\varphi) f(\lambda)\frac{d\lambda}{\lambda^4},\notag
	\end{equation}
		\begin{align}	
		K(R\lambda;p,q,\varphi)&=\lambda^4\int\limits_{S_m}\int\limits_{S_n} J_0(\rho\lambda) dS_mdS_n\notag\\	&=\lambda^4\int\limits_{x_n}^{x_n+h_x}\int\limits_{y_n}^{y_n+h_y}\int\limits_{x_m}^{x_m+h_x}\int\limits_{y_m}^{y_m+h_y}J_0(\rho\lambda)dx_0dy_0dxdy\notag\\
		&=\int\limits_{0}^{R\lambda p}\int\limits_{0}^{R\lambda q}\int\limits_{R\lambda\left(\cos\varphi-\frac{p}{2}\right)}^{R\lambda\left(\cos\varphi+\frac{p}{2}\right)}\int\limits_{R\lambda\left(\sin\varphi-\frac{q}{2}\right)}^{R\lambda\left(\sin\varphi+\frac{q}{2}\right)}J_0(\tau)d\tilde x d\tilde x_0 d\tilde y d\tilde y_0,\label{transf}
			\end{align}	where
	\begin{align}
	  &\tau=\sqrt{(\tilde x-\tilde x_0)^2+(\tilde y-\tilde y_0)^2};\notag\\
		&R=\sqrt{\left(x_n-x_m-\frac{h_x}{2}\right)^2+\left(y_n-y_m-\frac{h_y}{2}\right)^2};\notag\\
		&p=\frac{h_x}{R}~ ~ ~q=\frac{h_y}{R}~ ~ ~\varphi=\arctan\frac{y_n-y_m-\frac{h_y}{2}}{x_n-x_m-\frac{h_x}{2}}.\notag
	\end{align}
	Let $\lambda=e^{-t}$, $R=e^{s}$
	\begin{equation}
		\label{convolv_app}
		\begin{aligned}
		F\left(e^{s};p,q,\varphi\right)=e^{3s}\int\limits_{-\infty}^{\infty}\Phi(s-t;p,q,\varphi)f\left(e^{-t}\right)dt,\\
		\Phi(s-t;p,q,\varphi)=K\left(e^{s-t};p,q,\varphi\right)e^{-3(s-t)}.
	\end{aligned}
	\end{equation}
	For fixed $p,q,\varphi$ the integral in \eqref{convolv_app} is the convolution integral with kernel $\Phi$. Note that for different values of $\alpha$ and $\beta$ the kernels can be obtained similarly.

	 The main advantage of using the convolution integrals  is  that their computation does not require the explicit calculation of kernel $\Phi$. 
Consider the input function $v(t)$ and the output function $u(s)$ such that
	\begin{equation}
		u(s)=\int\limits_{-\infty}^{+\infty}\Phi(s-t)v(t)dt.
	\end{equation}
	For some  $N=2M$, $l$, $0<\xi<0.5l$, $k=-M\dots M-1$, define
	\begin{equation}
		\label{W}
		W_s=(-1)^s\frac{1}{N}\sum\limits_{n=-M}^{M-1}\left\{\frac{\sum\limits_{m=-M}^{M-1}(-1)^m u(ml-\xi)e^{-2i\pi \frac{mn}{N}} }{\sum\limits_{m=-M}^{M-1}(-1)^m v(ml+0.5l)e^{-2i\pi \frac{mn}{N} } }\right\}e^{2i\pi\frac{sn}{N}}.
	\end{equation}
	Then for any $g(t)$ one gets
	\begin{equation}
		\label{quad1}
		\int\limits_{-\infty}^{+\infty}\Phi(s-t)g(t)dt \approx \sum\limits_{s=N_1}^{N_2} W_s g(sl+\xi),~ ~ ~ -M\le N_1 < N_2 \le M-1.
	\end{equation}
	The tradeoff between the accuracy of \eqref{quad1} and the  computational time is achieved by the particular selection of  $M$, $N_1$, $N_2$, $l$, $\xi$ and functions $u$ and $v$.

	From~\eqref{transf}~and~\eqref{quad1} the approximation formulas~for~\eqref{g_eint} can be obtained
	\begin{equation}
		\label{main_q1}
		\begin{aligned}
			&I_{\alpha,\beta} \approx R^{(3-\alpha-\beta)} \sum\limits_{m=N_1}^{m=N_2}W_m^{\alpha,\beta}(p,q,\varphi,\alpha,\beta)f\left(\frac{\lambda_m}{R}\right),\\
		&\lambda_m=e^{ml+\xi}.
		\end{aligned}
	\end{equation}

	For the given  input function  $v(t)=8e^{-t^2}\left(t^5-4t^3+2t\right)$ the output functions for different kernels  can be expressed analytically by Gaussian and error functions.
	Inspired by \cite{Anderson1979} the parameters used are $l=0.2$, $\xi=0.0964$, $M=512$, $N_1=-250$, $N_2=200$. In computational experiments these parameters provided appropriate accuracy in calculation of  $\widehat B_n^m$.

\end{appendices} 

\end{document}